\documentclass[a4paper,twoside,absolute]{mymathart}
\usepackage{enumerate}
\usepackage[theorems_new]{mymathmacros}
\usepackage{subfigure}
\usepackage{ogonek}

\numberwithin{equation}{section}

\usepackage{hyperref}

\renewcommand{\H}{\mathbb{H}}

\renewcommand{\theta}{\vartheta}
\renewcommand{\phi}{\varphi}
\renewcommand{\rho}{\varrho}

\newcommand{\adds}{\underline{s}}

\newcommand{\s}{\adds}

\newcommand{\bdyit}[2]
             {{\rule{0pt}{0pt}_{\mbox{$\scriptstyle #2$}}^{\mbox{%
                   $\scriptstyle #1$}} }}

\newcommand{\W}{\mathcal{W}}

\renewcommand{\P}{\mathcal{P}}

  \newcommand{\B}{\mathcal{B}}
  \renewcommand{\S}{\mathcal{S}}

\newcommand{\Blog}{\B_{\log}}

\newcommand{\V}{\mathcal{V}}
\newcommand{\U}{\mathcal{U}}

\newcommand{\id}{\operatorname{id}}
\newcommand{\dil}{\operatorname{dil}}

\renewcommand{\picturedir}{.}

\renewcommand{\subproofsymbol}{}

\title[Rigidity of Escaping Dynamics]{%
  Rigidity of Escaping Dynamics \\
  for Transcendental Entire Functions}
\author{Lasse Rempe}
\address{Department of Mathematical Sciences, University of Liverpool, L69 7ZL,
United Kingdom}
\email{l.rempe@liverpool.ac.uk}
\date{\today}
\thanks{Supported in part
 by a postdoctoral fellowship of the 
 German Academic Exchange Service (DAAD) and by EPSRC
 Advanced Research Fellowship EP/E052851/1}

\subjclass[2000]{Primary 37F10; Secondary 30D05}

\begin{document}

\begin{abstract}
 We prove an analog of B\"ottcher's theorem for 
  transcendental entire functions in the Eremenko-Lyubich class
  $\B$. More
  precisely, let $f$ and $g$ be entire functions
  with bounded sets of singular values and suppose that $f$ and $g$
  belong to the same parameter space
  (i.e., are \emph{quasiconformally equivalent} in the sense
  of Eremenko and Lyubich). Then $f$ and $g$ are conjugate when
  restricted to the set of points that remain in some sufficiently
  small neighborhood of infinity under iteration. Furthermore,
  this conjugacy extends to a quasiconformal self-map of the plane.

 We also prove that the conjugacy is essentially unique. In particular,
  we show that a function $f\in\B$ has no invariant line fields on its
  escaping set.
  Finally, we show that any two hyperbolic functions $f,g\in\B$ that
  belong to the same parameter space are conjugate on their sets of
  escaping points. 
\end{abstract}

\maketitle

\section{Introduction}

 The study of the dynamical behavior of transcendental functions,
  initiated by Fatou in 1926 \cite{fatou}, has enjoyed increasing interest 
  recently. Many intriguing phenomena discovered in
  polynomial dynamics, relating to the
  behavior of high-order renormalizations of a polynomial, 
  occur naturally for transcendental maps. Compare, for example,
  Shishikura's proof that the boundary of the Mandelbrot set has
  Hausdorff-dimension $2$ \cite{mitsudim} with McMullen's 
  treatment of the Julia set of $z\mapsto\lambda\exp(z)$ 
  \cite{hausdorffmcmullen}. A
  more recent example is provided by work of
  Avila and Lyubich \cite{avilalyubichfeigenbaum2}, who proved that
  a
  constant-type Feigenbaum quadratic 
  polynomial with positive measure Julia set would have hyperbolic 
  dimension less than $2$. Work of Urbanski and Zdunik \cite{urbanskizdunik1} 
  shows that a similar
  phenomenon occurs for the simplest exponential maps. 

\begin{figure}%
  \subfigure[$f_1(z)= 2(\exp(z)-1)$]{\hspace{.5cm}%
   \includegraphics[height=.35\textheight]{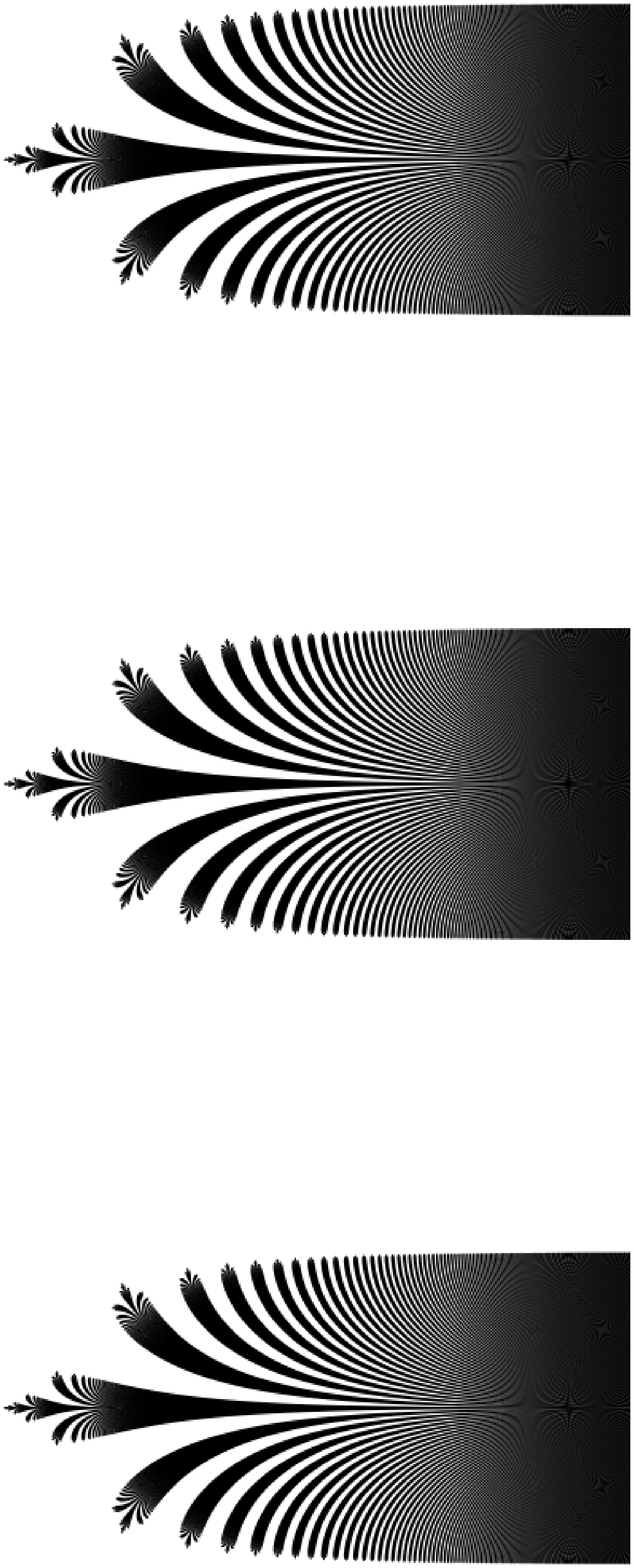}\hspace{.5cm}}
    \hfill%
  \subfigure[$f_2(z)= (z+1)\exp(z)-1$]{\hspace{.5cm}
    \includegraphics[height=.35\textheight]{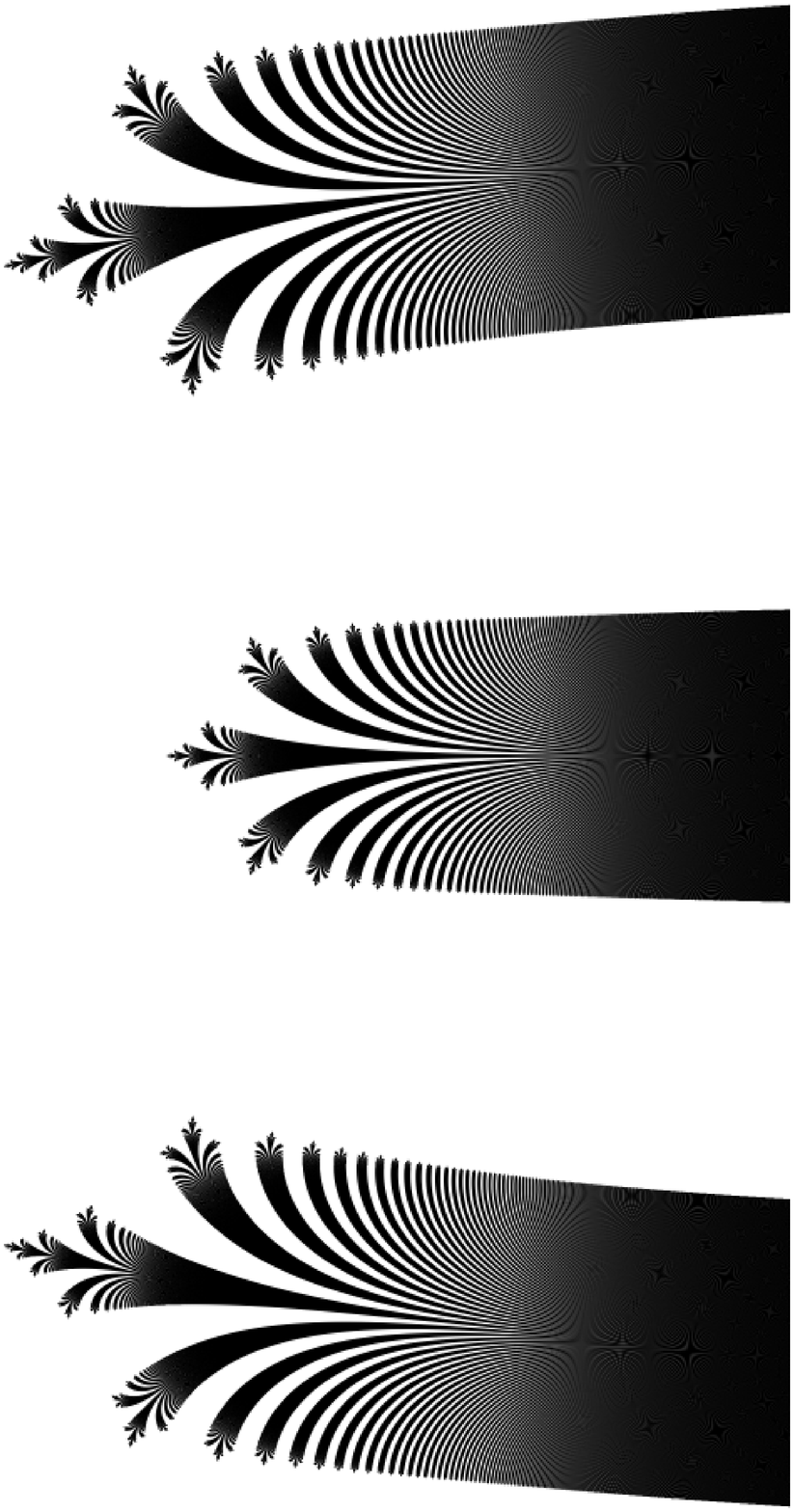}\label{fig:zexp}\hspace{.5cm}}%
 \hfill%
  \subfigure[$f_3(z)= \lambda\sinh(z)$]{\hspace{.5cm}
 \includegraphics[height=.35\textheight]{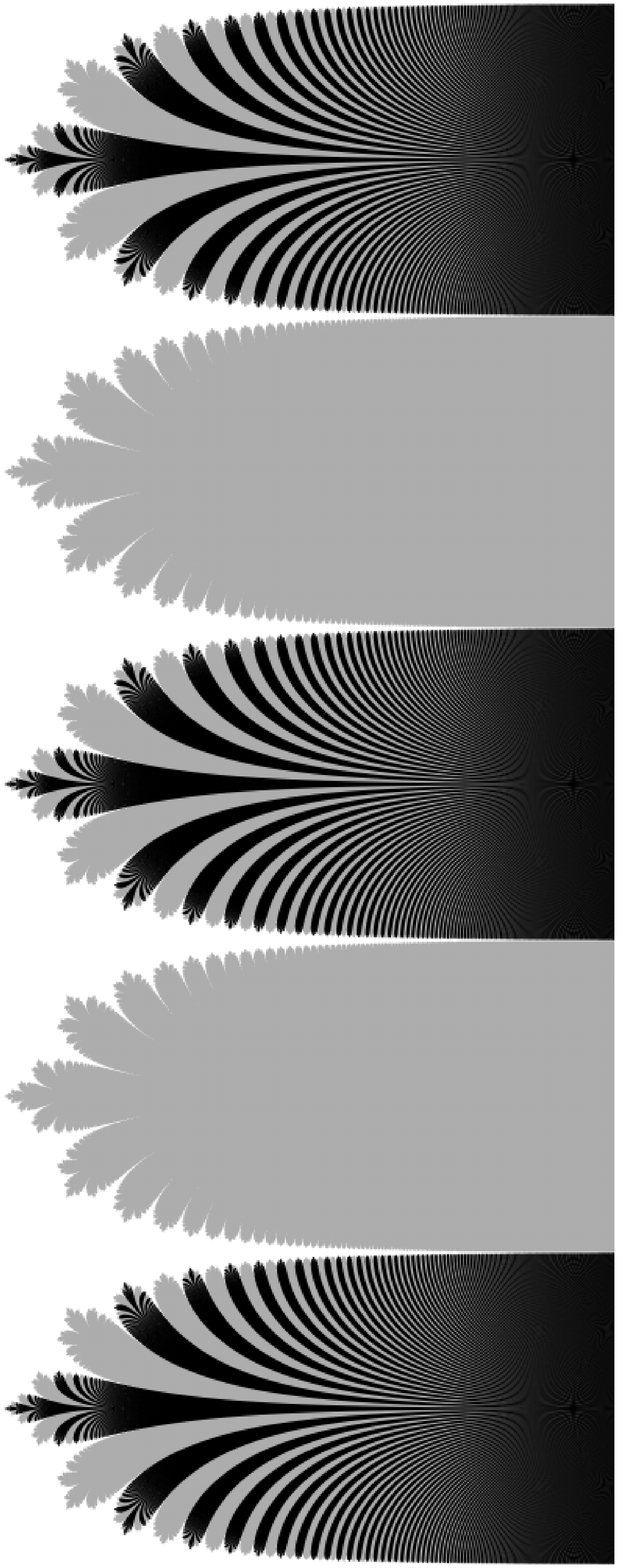}\hspace{.5cm}%
               \label{fig:sinhexp}%
           \hspace{.5cm}}
  \caption{\label{fig:julia} Images (a) and (b) show the Julia sets
   of the functions $f_1$ and $f_2$ (in black). Our results imply that
   these two functions are quasiconformally conjugate in a neighborhood of
   these sets. (Compare Theorem \ref{thm:disjointconjugacy} and
   Remark \ref{rmk:quasidisks}.) In (c) the black set consists of 
   points whose orbits under $f_3$ remain in a right half
   plane. Again, restricted to this set, $f_3$ is quasiconformally conjugate
   to $f_1$. (The Julia set of $f_3$ is underlaid in gray.) 
   Note that the three maps are function-theoretically diverse:
   $f_1$ has one asymptotic value, 
   $f_2$ has both an asymptotic and a critical value, and
   $f_3$ has two critical values. 
   (In (c), $\lambda=0.575$.)}
 \end{figure}

 In this note, we prove a structural theorem for
  the dynamics near a logarithmic singularity. On the one
  hand, this result explains the
  observation that many Julia sets of explicit 
  entire transcendental functions 
  bear striking similarities to each
  other, even if they are very different from a
  function-theoretic point of view, compare Figure \ref{fig:julia}. 
  On the other hand,
  it provides a tool to better understand the Julia sets of
  these functions, and results in some important rigidity
  statements required in the study of density of hyperbolicity
  \cite{lassesebastiandensity}. 

 The \emph{Eremenko-Lyubich class} $\B$ is the class of
  transcendental entire functions for which the set
  $\sing(f^{-1})$ of critical and asymptotic values is bounded. 
  We say that two functions $f,g\in\B$ are
  \emph{quasiconformally equivalent near $\infty$} if
  there exist quasiconformal maps $\phi,\psi:\C\to\C$ such that
   \begin{equation}
    \psi(f(z)) = g(\phi(z))  \label{eqn:equivalence}
   \end{equation}
  whenever $|f(z)|$ or $|g(\phi(z))|$ is large enough.
  (When (\ref{eqn:equivalence}) holds on all of 
   $\C$, the maps are called \emph{quasiconformally equivalent};
   compare Section
   \ref{sec:preliminaries}. Quasiconformal equivalence
   classes form the natural parameter spaces of entire functions.)

 \begin{thm}[Conjugacy near infinity] \label{thm:main1}
  Let $f,g\in\B$ be quasiconformally
   equivalent near infinity. Then there exist $R>0$ and
   a quasiconformal map
   $\theta:\C\to\C$ such that
    $\theta\circ f = g\circ \theta$
    on 
   \[J_R(f) :=
      \{z:|f^n(z)|\geq R \text{ for all $n\geq 1$}\}. \]
   Furthermore, 
   $\theta$ has
   zero dilatation on $\{z\in J_R(f):|f^n(z)|\to\infty\}$.
 \end{thm}
 \begin{remark}[Remark 1]
   In fact, our methods are 
    purely local and as such apply to
    any (not necessarily globally defined) function that has
    only logarithmic singularities over infinity. In particular,
    they apply to restrictions of certain entire (or meromorphic)
    functions that themselves do not belong to class $\B$.     
    We refer the reader to Section \ref{sec:preliminaries} for
    the precise definition of the class of functions that is treated.
 \end{remark}
 \begin{remark}[Remark 2]
  For functions with non-logarithmic singularities over infinity,
    the dynamics near infinity may vary dramatically within the
    same parameter space. For example, for the function 
    $z\mapsto z-1-\exp(z)$, all points with sufficiently negative real part
    tend to $-\infty$ under iteration: the function has
    a \emph{Baker domain} containing a left half plane. On the
    other hand, the function $z\mapsto z+1-\exp(z)$ does not have any 
    Baker domains: every orbit in the Fatou set converges to
    an attracting fixed point; see \cite[Section 5.3]{weinreichthesis}.
 \end{remark}

 \begin{figure}
  \subfigure[$f_4(z) = \exp(z)+\kappa$]{%
\includegraphics[width=.45\textwidth]{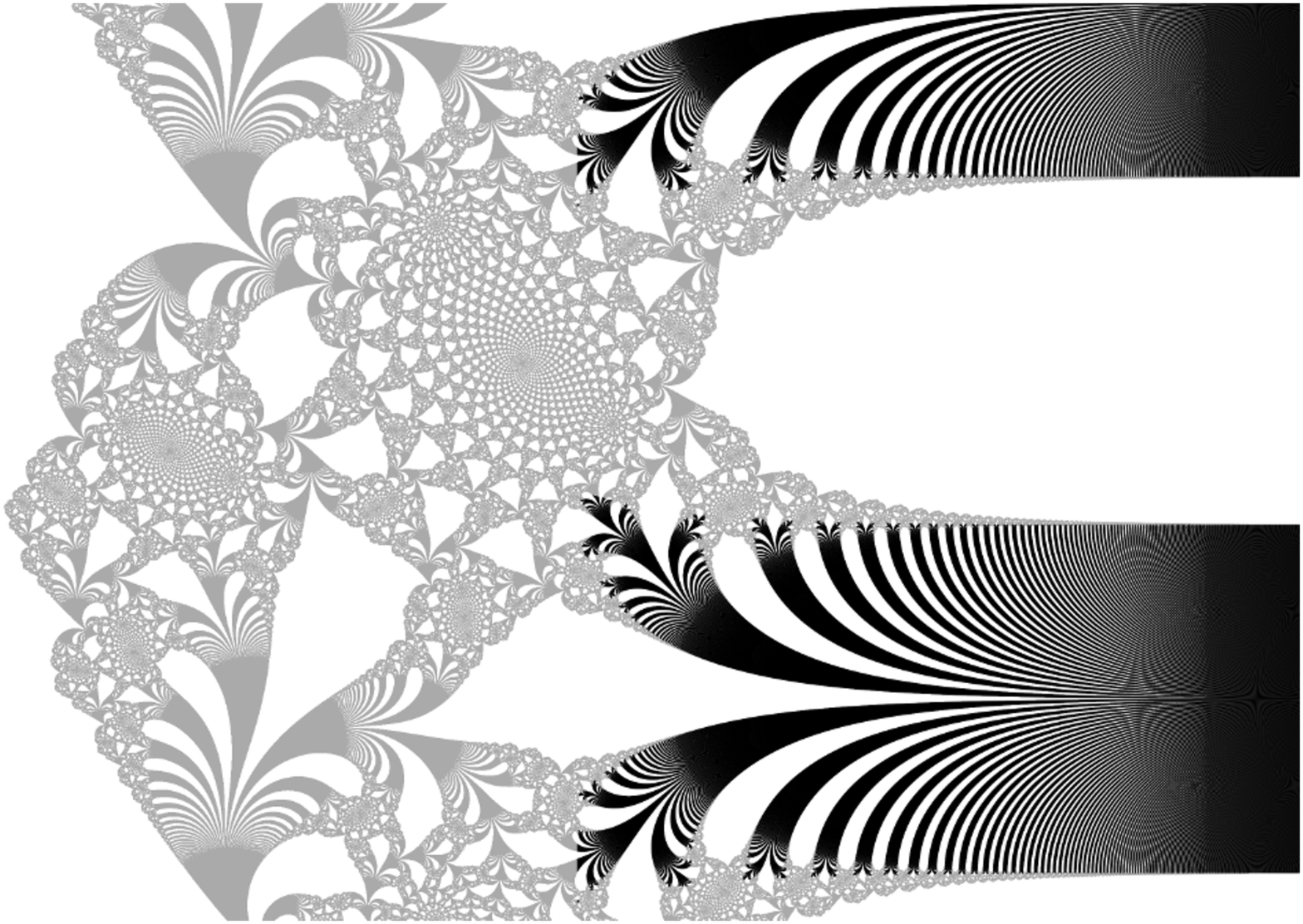}\label{fig:exp_attracting}}%
   \hfill
  \subfigure[$f_5(z)= \exp(z)$]{%
\includegraphics[width=.45\textwidth]{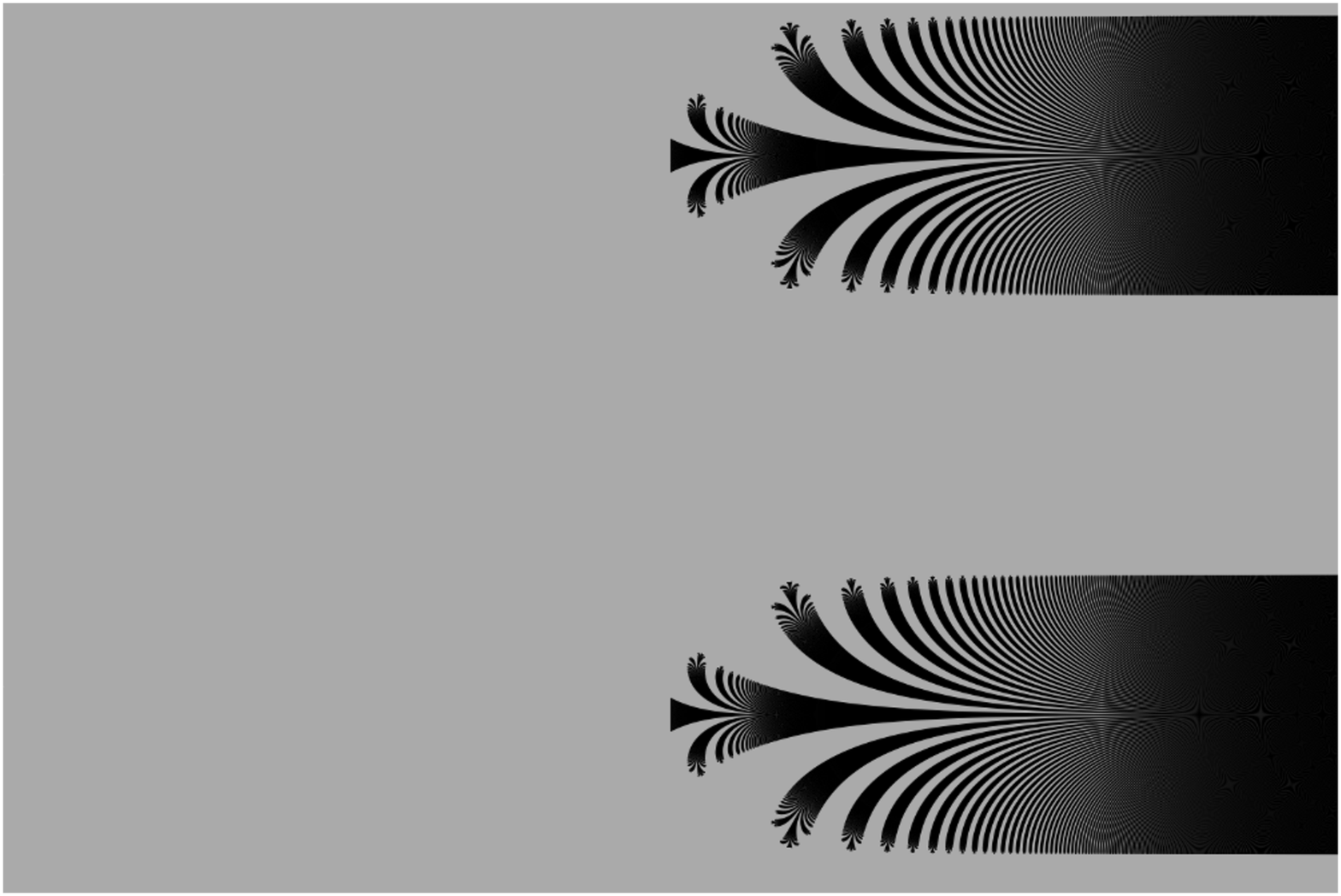}%
   \label{fig:sinhz_z}}%
 \caption{\label{fig:JR}Two functions that are quasiconformally 
   equivalent to the map
   $f_1$ from figure \ref{fig:julia}, but have very different dynamics:
   in (a), the Julia set (in gray) is a ``pinched Cantor bouquet'', while
   in (b) it is the entire complex plane. However, on the sets
   $J_R(f_i)$ from Theorem \ref{thm:main1} (in black), they are 
   quasiconformally conjugate to (a suitable restriction of) $f_1$.
   (The parameter in (a) is given by $\kappa=1.0038+2.8999i$.)}
 \end{figure}

 Theorem \ref{thm:main1} can be seen as an analog to a classical
  theorem of B\"ottcher  (see \cite[Theorem 18.10]{jackdynamicsthird}),
  stating that any two polynomials of the same degree $d\geq 2$
  are conformally conjugate near $\infty$. 
  We find the generality of our theorem surprising for a number of
  reasons. Not only can functions that are quasiconformally equivalent
  near infinity have very different function-theoretic properties
  (recall Figure \ref{fig:julia}), but more significantly the behavior 
  near infinity can vary widely between different functions in $\B$. 
  Indeed, for the function-theoretically
  simplest functions in this class, such as those shown in
  Figure \ref{fig:julia},  and in fact all functions $f\in\B$
  of finite order \cite[Theorem 1.2]{strahlen},
  the escaping set 
    \begin{equation}
      I(f) := \{z\in\C: f^n(z)\to\infty \} \label{eqn:escapingset}
    \end{equation}
  consists entirely of curves. On the other hand,
  it is is possible for the escaping set of a hyperbolic
  function $f\in \B$ to contain no nontrivial curves at all
   \cite[Theorem 8.4]{strahlen}. 
  Theorem \ref{thm:main1} shows that, even for such
    a ``pathological'' function, 
    the behavior near infinity remains
    the same throughout its quasiconformal equivalence class. 

 Douady and 
  Hubbard \cite{orsay} used B\"ottcher's theorem to introduce
  \emph{dynamic rays}, 
  which have become the backbone of the
  successful theory of polynomial dynamics. We believe that
  our result will likewise be useful in the study of families
  of transcendental functions, even those with such wild behavior as 
  the example mentioned above. Indeed, one corollary of Theorem 
  \ref{thm:main1}
  (Corollary \ref{cor:nocurves}) is that any function that is
  quasiconformally equivalent to this example also contains no
  curves in its escaping set. 

\smallskip

 Another aspect of the theorem's generality
  that seems surprising
  is the statement about dilatation. It is worth noting that two 
  quasiconformally
  equivalent functions in class $\B$ may have different orders of growth.
  (Whether this is possible for functions with \emph{finitely many} singular
   values is a difficult open problem.) Hence the map $\phi$ in
   the definition of quasiconformal equivalence cannot, in general, be chosen
   to be asymptotically conformal near infinity. In such a situation, one
   could imagine that some of the dilatation of the quasiconformal map $\theta$
   would 
   be supported on the escaping set $I(f)$, 
   but by Theorem \ref{thm:main1} this is not the case.

  In fact, we will show that the map $\theta$ is
   essentially unique (more precisely, it is unique up to 
   an initial choice of isotopy class; compare Corollary \ref{cor:uniqueness});
    hence it follows that \emph{no} quasiconformal
   conjugacy between $f$ and $g$ can support dilatation on the set $I(f)$.
 \begin{thm}[No invariant line fields] \label{thm:nolinefields}
  A function $f\in\B$ supports no invariant line fields on its
   escaping set.
 \end{thm}
 \begin{remark}[Remark 1]
  This statement has content only in families 
   where 
   the set of escaping points has positive measure. As far as we know, 
   it is new even for the family $z\mapsto a\exp(z) + b\exp(-z)$ of cosine
   maps, whose escaping sets have positive measure by \cite{hausdorffmcmullen}.
 \end{remark}
 \begin{remark}[Remark 2]
  Showing that the Julia set of a polynomial cannot support
     an invariant line field is a major open problem in complex dynamics.
     In contrast, it is known \cite{alexmishaexamples} that
     there are entire functions with invariant line fields on their Julia sets.
     In fact, the example from \cite{alexmishaexamples} has an invariant
     line field on $I(f)\cap J(f)$, showing that Theorem \ref{thm:nolinefields}
     becomes false if the assumption $f\in\B$ is dropped. 
 \end{remark} 

 By the same reasoning, we also obtain further 
   rigidity principles for the set $I(f)$, of which the following
   is an important special case. 

 \begin{thm}[QC rigidity on escaping orbits]
   \label{thm:qcrigidity}
  Suppose that $f$ and $g$ are entire functions with finitely many
   singular values, and let $\pi$ be a topological
   conjugacy betweeen $f$ and $g$. If 
   $\mathcal{O} = \{z_0, f(z_0),f^2(z_0),\dots\}$
   is any escaping orbit of $f$, then the restriction
   $\pi|_{\mathcal{O}}$ extends to a quasiconformal self-map of the plane. 
 \end{thm}

 While the source of the
   rigidity here is
   much softer than in the famous rigidity results
   for rational functions (as indicated by the absence of dynamical
   hypotheses), our results provide an essential step in 
   transferring 
   rigidity theorems from the rational 
   to the transcendental setting. For example, in 
   \cite{linefields}, Theorem \ref{thm:nolinefields} is used to obtain
   the 
   absence of invariant line fields on the Julia sets of 
   a large class of ``nonrecurrent'' transcendental functions,
   extending work of
   Graczyk, Kotus and {\'S}wi\k{a}tek \cite{graczykkotusswiatek}. In
   \cite{lassesebastiandensity}, our results are used, together with
   work of Kozlovski, Shen and van Strien 
    \cite{densityaxioma,densityaxiomadimensionone} to establish
   density of hyperbolicity in certain families of real transcendental
   entire functions (including the real cosine family
   $a\sin(z)+b\cos(z)$, $a,b\in\R$).

 \smallskip

 In contrast to the polynomial case, the map $\theta$
  from Theorem \ref{thm:main1} will generally \emph{not}
  extend to a conjugacy between the escaping sets of 
  $f$ and $g$ \cite[Proposition 2.1]{topescapingnew}. However, in the case of
  \emph{hyperbolic} 
  functions $f\in\B$ --- i.e., those for which the postsingular set is
   compactly contained in the Fatou set ---  we can do better. 

 \begin{thm}[Conjugacy for hyperbolic maps]
   \label{thm:hyperbolicconjugacy}
  Let $f,g\in\B$ be quasiconformally equivalent near infinity, 
   and suppose that $f$ and $g$ are hyperbolic. 

  Then $f$ and $g$ are conjugate on their
   sets of escaping points.
 \end{thm}

 Together with recent results of Bara{\'n}ski \cite{baranskihyperbolic},
  our proof of Theorem \ref{thm:hyperbolicconjugacy}  also shows  
  that, for hyperbolic $f\in\B$ of finite order, $J(f)$ can be described
  as a \emph{pinched Cantor Bouquet}; i.e., as the quotient of
  a Cantor Bouquet (or ``straight brush'')
  by a closed equivalence relation on the endpoints. Recently,
  Mihaljevi\'c-Brandt \cite{helenaconjugacy}
  has generalized Theorem \ref{thm:hyperbolicconjugacy}
  to a large class of ``subhyperbolic'' entire functions. In particular,
  her result applies to all postcritically finite 
  functions $f\in\B$ with no asymptotic values for which there is some 
  $\Delta$ such that all critical points of $f$ have degree at most $\Delta$.

 \subsection*{Structure of the article and ideas of the proofs}
  We begin in Section \ref{sec:preliminaries} by reviewing some basic
   properties of Eremenko-Lyubich functions and introducing the local 
   class
   $\Blog$. Section 
   \ref{sec:conjugacy} is devoted to the proof of Theorem
   \ref{thm:main1}, which has two main ingredients. 
    The first of these is 
    the well-known fact that functions in $\B$ are expanding inside
    their logarithmic tracts. The second is that
    the quasiconformal maps
    $\phi$ and $\psi$ do not move points near infinity
    more than a finite distance
    with respect to the hyperbolic metric in a punctured neighborhood of 
    infinity. With these two facts, most of
    the theorem can be considered to be a variant
    of standard conjugacy results for expanding maps. 

  However, in
    order to obtain the statement on dilatation, we need to break the
    proof down into two cases: one where both maps $f$ and $g$ are
    dynamically simple (``disjoint-type'') functions, and one where the 
    quasiconformal maps $\phi$ and $\psi$ are in fact affine. (In the 
    latter case, the quasiconformality of the function $\theta$, and the
    dilatation estimate, will be obtained via the ``$\lambda$-lemma'' of
    \cite{mss}.)
    Together, these two cases combine to give the full theorem;
    compare also the discussion at the end of Section \ref{sec:conjugacy}.

   The proofs of Theorems \ref{thm:nolinefields} and \ref{thm:qcrigidity}
    are given in Section \ref{sec:rigidity}. As already mentioned, they
    rely on the fact that the map $\theta$ is unique in a certain sense
    (Corollary \ref{cor:uniqueness}). 
    The idea of the proof can be 
    traced back to the argument of Douady and Goldberg \cite{douadygoldberg}
    who proved that two topologically conjugate
    real exponential maps with escaping singular orbits
    must be conformally conjugate. 
  
  To prove Theorem \ref{thm:hyperbolicconjugacy}
   in Section \ref{sec:hyperbolic}, we show that hyperbolic
   entire functions are expanding with respect to the hyperbolic
   metric; the construction of a semi-conjugacy then proceeds as usual
   for expanding maps.

  In Appendix \ref{sec:fatoueremenko}, we discuss the relation of our results
   with some well-known questions regarding escaping sets posed by
   Fatou \cite{fatou} and Eremenko \cite{alexescaping}. 

 \subsection*{Acknowledgements} 
  I would like to especially thank
   Carsten Petersen, whose
   thought-provoking questions during
   a talk at the Institut Henri Poincar\'e
   initiated the research that led
   to these results.
   I would also 
   like to thank Walter Bergweiler, Adam Epstein, 
   Alex Eremenko, Jeremy Kahn, Misha Lyubich,
   Phil Rippon, 
   Dierk Schleicher, Gwyneth Stallard,
   Sebastian van Strien and particularly Helena Mihaljevi\'c-Brandt for 
   interesting discussions and comments, and the referee for helpful 
   suggestions 
   that led to marked improvement in exposition. 
  
 \subsection*{Background and Notation}
  We refer the reader to 
   \cite{jackdynamicsthird,waltermero,hubbardteichmueller,lehtovirtanen} for
   introductions
   to holomorphic dynamics, plane 
   hyperbolic geometry and
   quasiconformal mappings. 

  We denote the complex plane by 
   $\C$ and the Riemann sphere by
   $\Ch = \C\cup\{\infty\}$. All closures and boundaries will be
   understood to be taken in $\C$, unless explicitly stated otherwise. 
   We denote the right 
    half-plane by
    $\H := \{\re z >0\}$;  
   more generally, we write
   \[ \H_{Q} := \{\re z > Q\}.
   \]

  If $f:\C\to\C$ is an entire function, we denote its
   \emph{Julia} and \emph{Fatou} sets by $J(f)$ and $F(f)$,
   respectively. Recall that the \emph{escaping set} $I(f)$ was defined in
   (\ref{eqn:escapingset}). 

  The set of \emph{singular values}, $S(f)$, is the closure of the
   set $\sing(f^{-1})$ of critical and asymptotic values of $f$. 
   The Speiser class $\S$ and the Eremenko-Lyubich class
   $\B\supset \S$ are defined as
   \begin{align*}
      \S &:= \{f:\C\to\C\ \text{entire, transcendental}:
                 S(f)\text{ finite} \} \quad\text{and}\\
      \B &:= \{f:\C\to\C\ \text{entire, transcendental}:
                  S(f)\text{ bounded}\}.
   \end{align*}

\section{Preliminaries} \label{sec:preliminaries}

 \subsection*{The hyperbolic metric}
  If $U\subset\C$ is open and $\C\setminus U$ contains 
  at least two points, we denote the
  density of the hyperbolic metric in $U$ by 
  $\rho_U$.  We denote hyperbolic distance and length in $U$ by
  $\dist_U$ and $\ell_U$, respectively. The derivative of a 
  holomorphic function $f$
  with respect to the hyperbolic metric of $U$ (where defined)
   will be denoted by
   \[ \|Df(z)\|_U := |f'(z)| \cdot \frac{\rho_U(f(z))}{\rho_U(z)}. \]

 Recall \cite[Corollary A.8]{jackdynamicsthird} 
   that, if $U$ is simply connected, then
    \begin{equation}
      \frac{1}{2\dist(z,\partial U)} \leq
         \rho_U(z) \leq \frac{2}{\dist(z,\partial U)} 
              \label{eqn:standardestimate}
    \end{equation}
  for all $z \in U$; we refer to this
  as the \emph{standard estimate} on the hyperbolic metric. We also remind the
  reader that holomorphic covering maps preserve the hyperbolic metric,
  and that Pick's theorem \cite[Theorem 2.11]{jackdynamicsthird} states that
  $\rho_{U'}(z) > \rho_{U}(z)$ for all $z\in U'$ if $U'\subsetneq U$. 

 In Section \ref{sec:hyperbolic}, we will use the following
  estimate on the hyperbolic metric in certain multiply-connected domains.
 \begin{lem}[Hyperbolic metric in countably punctured sphere]
    \label{lem:beardonpommerenke}
  Let $(w_j)_{j\in\N}$ be a sequence of points in $\C\setminus\{0\}$ with
    $w_j\to\infty$ and satisfying
     $|w_{j+1}|\leq C|w_j|$
     for some constant $C>1$ and all $j\in\N$. Set
     $V := \C\setminus (\{0\}\cup\{w_j:j\in\N\})$. Then
     $1/\rho_{V}(z) = O(|z|)$ as $z\to\infty$. 
 \end{lem}
 \begin{proof}[Sketch of proof]
  We use the following estimate on the hyperbolic metric in a
   doubly punctured plane $U_{a,b} := \C\setminus\{a,b\}$: if 
   $|z-a|\leq |z-b|$, then 
     \begin{equation} \label{eqn:beardonpommerenke}
       1/\rho_{U_{a,b}}(z) \leq
          K\cdot |z-a|\cdot \left(1+\left|\log\frac{|b-a|}{|z-a|}\right|\right). \end{equation}
   (In 
   \cite{beardonpommerenke}, this expression is used to
   describe the precise order of magnitude for
   the hyperbolic density of an arbitrary multiply-connected domain.)

  To prove the claim, let $z\in V$ with $|z|\geq |w_0|$
    and let $a\in\partial V$ such that $|z-a|$ is minimal.
    We obtain the desired statement by applying (\ref{eqn:beardonpommerenke})
    with a suitable point $b\in\partial V$ 
    as described below, and using Pick's theorem. 
   \begin{itemize}
    \item If $a = 0$, let $b=w_j$ where $j$ is minimal with $|w_j|\geq |z|$;
       we obtain
       \[ 1/\rho_V(z) \leq 1/\rho_{U_{a,b}}(z) \leq
             K\cdot |z|\cdot (1+\log C). \]
    \item If $a\neq 0$ and $|z-a| > |a|/2$, let $b=w_j$, where $j$ 
      is minimal with $|w_j|\geq 3|z-a|$. Then 
       \[ 1/\rho_V(z) \leq 1/\rho_{U_{a,b}}(z) \leq
             K\cdot |z|\cdot (1+ \log(3C)). \]
    \item If $a\neq 0$ and $|z-a|\leq |a|/2$, let $b=0$. In this case,
       \begin{align*}
         1/\rho_V(z) & \leq 1/\rho_{U_{a,b}}(z) \leq
             K\cdot |z-a|\cdot (1+\log(|a|/|z-a|)) \\
          &\leq K \cdot |z-a| \cdot |a|/|z-a| =
             K\cdot |a| \leq 2K\cdot |z|. \end{align*}    
   \end{itemize}
   So overall we have $1/\rho_V(z) \leq (1+\log(3C))K|z| = O(|z|)$, 
     as claimed. 
 \end{proof} 

 \subsection*{Tracts and logarithmic coordinates}
  A domain $U\subset\C$ is called an \emph{unbounded Jordan domain} if
   the boundary of $U$ on the Riemann sphere is a Jordan curve passing
   through $\infty$. 

  Suppose that $f\in\B$, and let $D\subset \C$ be a bounded
   Jordan domain chosen such that 
   $S(f)\cup\{0,f(0)\}\subset D$. (E.g., 
   $D=\D_R(0)$, where $R \geq 1+|f(0)| + 
     \max_{s\in S(f)}|s|$.)  
   Let us set $W := \C\setminus \cl{D}$ and
   $\U := f^{-1}(W)$. Then
   each component $T$ of $\U$ is 
   an unbounded
   Jordan domain
   (called a \emph{tract} of
    $f$), and
  $f:T\to W$ is a universal covering.

   We can perform a \emph{logarithmic change of
    coordinates} (see \cite[Section 2]{alexmisha} or
     \cite[Section 4.8]{waltermero}) to obtain a $2\pi i$-periodic function 
     $F:\V\to H$,
    where $H=\exp^{-1}(W)$ and $\V=\exp^{-1}(\U)$, 
    such that $\exp\circ F = f\circ \exp$. We will say that this function
    $F$ is a \emph{logarithmic transform} of $f$.  By construction, 
    the following properties hold.

   \begin{enumerate}[(A)] 
    \item $H$ is a $2\pi i$-periodic unbounded Jordan domain that contains
     a right half-plane.   \label{item:H}
    \item $\V\neq\emptyset$ 
     is $2\pi i$-periodic and $\re z$ is bounded from below in
     $\V$.
    \item $F$ is $2\pi i$-periodic. \label{item:periodicity}
    \item Each component $T$ of $\V$ is an 
      unbounded Jordan domain that is disjoint
      from all its $2\pi i\Z$-translates. 
      For each such $T$, the restriction
      $F:T\to H$ is a conformal 
      isomorphism with $F(\infty)=\infty$. 
      ($T$ is called a \emph{tract of $F$}; we 
       denote the inverse of $F|_T$ by 
       $F_T^{-1}$.) \label{item:tracts}
    \item The components of $\V$ accumulate only at $\infty$; i.e.,
      if $z_n\in\V$ is a sequence of points all belonging to different
      components of $\V$, then $z_n\to\infty$.  \label{item:accumulatingatinfty}   \end{enumerate}

   We will denote by $\Blog$ the class of all functions
    \[ F:\V\to H, \]
   where $F$, $\V$ and $H$ have the properties
   (\ref{item:H}) to (\ref{item:accumulatingatinfty}), regardless of whether
    $F$ arises as the logarithmic transform of a function $f\in\B$ or not.    
   
   Note that
    any
    $F\in\Blog$ extends continuously to $\cl{\V}$ by Carath\'eodory's theorem.
    The
    \emph{Julia set} and \emph{escaping set} of 
    $F\in\Blog$ are defined to be
   \begin{align*}
     J(F) &:= \{z\in \cl{\V}: F^n(z)\in\cl{\V}\text{ for all $n\geq 0$}\}
     \quad\text{and} \\
     I(F) &:= \{z\in J(F): \re F^n(z)\to\infty\}.
   \end{align*}
    When $F$ is the logarithmic transform of
    a function $f\in\B$, then
    $\exp(I(F))\subset I(f)$ and the orbit of
    every $z\in I(f)$ will eventually remain in $\exp(I(F))$. 
    For $Q>0$, we also define 
   \begin{align*} J_{Q}(F)& :=
     \{z\in J(F):
       \re F^n(z) \geq Q
       \text{ for all $n\geq 1$}\}
    \quad\text{and}\\
     I_{Q}(F) &:= I(F)\cap J_{Q}(F). 
    \end{align*}
  If $F$ is the logarithmic transform of $f$, then 
   clearly $\exp(J_Q(F))=J_{e^Q}(f)$ (the latter set was defined
   in Theorem \ref{thm:main1}). 

\subsection*{Expansion and normalization}
  Let us introduce two important sub-classes of $\Blog$.

  \begin{defn}[Disjoint-type and normalized functions]
   Let $F:\V\to H$ belong to the class $\Blog$.
  \begin{enumerate}
   \item
    We say that $F$ is of \emph{disjoint type} if
      $\cl{\V}\subset H$. 
  \item We say that $F$ is \emph{normalized} if
     $H=\H$  and 
     \begin{equation} \label{eqn:expansion}
      |F'(z)|\geq 2
     \end{equation} 
     for all $z\in\V$. 
   \end{enumerate}
  \end{defn}
  \begin{remark}
   If an entire function $f\in\B$ 
    has a logarithmic transform $F$ of disjoint type, then
    we will also say that $f$ itself is of disjoint type. In this case,
    the Fatou set of $f$ consists of a single immediate basin of
    attraction, and $J(f) = \exp(J(F))$. 
    The
    examples from Figure \ref{fig:julia} are of disjoint type, while
    those in Figure \ref{fig:JR} are not. 
  \end{remark}

 Let $F:\V\to H$ be any element of $\Blog$. It follows easily from
   (\ref{item:tracts}) and the standard estimate (\ref{eqn:standardestimate})
   on the hyperbolic metric that
    \begin{equation}
       \|DF(z)\|_H \to \infty \quad\text{as $\re(z) \to\infty$}. 
      \label{eqn:hypexpansion}
    \end{equation}
  In particular, by Pick's theorem, 
   any \emph{disjoint-type} function $F\in\B$ is uniformly expanding
   with respect to the hyperbolic metric in $H$. 

  The same argument also shows, again for any function $F\in\Blog$, that
   $|F'(z)|\to\infty$ as $\re(F(z))\to\infty$; see \cite[Lemma 1]{alexmisha}. 
   In particular, 
   there is $R>0$ such that (\ref{eqn:expansion}) holds 
   for all $z\in\V$ with $\re F(z) \geq R$. 
   By restricting $F$ to the set 
   $\wt{\V} := \{z\in\V: \re F(z) > R\}$ and 
   conjugating by $z\mapsto z-R$, we obtain the function 
     \[ \wt{F}: (\wt{\V}-R)\to \H;\quad
          z \mapsto F(z+R)-R. \]
   By construction, this function
    $\wt{F}$ is a normalized element of $\Blog$.
    As we are mostly concerned with the behavior of
    $F$ near $\infty$, we usually deal only with normalized
    functions. However, note that
    a normalization of a disjoint-type map $F$ need not be 
    of disjoint type.

  \begin{lem}[$J(F)$ has empty interior] \label{lem:emptyinterior}
   If $F\in\Blog$ is normalized or of disjoint type, then 
    $J(F)$ has empty interior.
  \end{lem}
  \begin{proof}[Sketch of proof]
   This is the same argument as in \cite[Theorem 1]{alexmisha},
   using the uniform expansion of the function $F$ in the
   Euclidean metric (in the normalized case) resp.\ the
   hyperbolic metric (for disjoint-type maps).
  \end{proof}
  \begin{remark}
   It follows that, for \emph{any} $F\in\Blog$,
    $J_Q(F)$ has empty interior for sufficiently large $Q$; if
    $F$ is the logarithmic transform of a function $f\in\B$, then
    similarly $\exp(J_Q(F))\subset J(f)$ for sufficiently large $Q$. 
  \end{remark}

 It is easy to see that $J_Q(F)\neq \emptyset$ for all $Q$; in fact,
  the following is true.

 \begin{prop}[Unbounded sets of escaping points 
  {\cite[Theorem 2.4]{landingsiegel}}] \label{prop:unboundedsets}
  Let $F\in\Blog$, and let $T$ be a tract of $F$. Then there is
   an unbounded, closed, connected set $A\subset T\cap I(F)$ such that 
     $\re F^j(z) \underset{j\to\infty}{\to} +\infty$
   uniformly on $A$. 
 \end{prop}
  \begin{remark}
   In \cite{landingsiegel}, the theorem is stated for entire functions in the
    Eremenko-Lyubich class, but the proof applies also to functions
    in $\Blog$. 
    It follows from the results of 
    \cite{eremenkoproperty} that  
    the set $A$ can be chosen to be forward-invariant, but we 
    do not require this.
    Compare \cite{walterphilgwyneth} for the existence of 
    unbounded connected sets of escaping points in more general situations.
  \end{remark}

 \begin{cor}[Density of escaping sets] \label{cor:unboundedescapingsets}
  Let $F\in\Blog$ and $Q\geq 0$. Then $I_Q(F)$ is nonempty, and $\re z$ is
   unbounded from above in $I_Q(F)$. 

  Furthermore,
   if $Q'\geq Q$ is sufficiently large, then
      \[ J_{Q'}(F) \subset \cl{I_Q(F)}. \]
 \end{cor}
 \begin{proof}[Sketch of proof]
  We may assume, once again, that $F$ is normalized. 
  The previous proposition implies
   that there is $Q'\geq Q + \pi/2$ such that, for every $M\geq Q'$, 
   there is a point $z\in I_{Q}(F)$ with $\re(z) = M$. 

  Let $z \in J_{Q'}(F)$, and note that $I_Q(F)$ is $2\pi i$-invariant.
   Therefore, for every $n\geq 1$ we can find 
   $w^n\in I_Q(F)$ with 
   $\re(w^n)=\re F^n(z)$ and $|\im(w^n) - \im F^n(z)|\leq \pi$. 

  Pulling $w^n$ back along the orbit of $z$, and using the
   expansion property (\ref{eqn:expansion}), we obtain a sequence of
   points $\omega^n \in I_Q(F)$ with $|\omega^n - z| \leq \frac{\pi}{2^n}$.
   Hence $z\in \cl{I_Q(F)}$, as required. 
 \end{proof}
   
 \subsection*{Quasiconformal equivalence}

 Following \cite[Section 3]{alexmisha},
 two entire
  functions $f,g\in\B$ are called \emph{quasiconformally equivalent}
  if there exist quasiconformal maps $\phi,\psi:\C\to\C$ with
   \begin{equation} \label{eqn:qcequivalence}
    g\circ\phi = \psi \circ f.
   \end{equation}
  The set
   of all functions $g$ that are quasiconformally equivalent
   to $f$ can be considered the natural parameter space of $f$. (If
   $S(f)$ is finite, then this set forms a finite-dimensional 
   complex manifold \cite[Section 3]{alexmisha}.)

 Similarly, let us say that two functions
  $F,G\in\Blog$ (with domains
  $\V$ and $\W$) are quasiconformally equivalent if there are
  quasiconformal maps $\Phi,\Psi:\C\to\C$ such that
  \begin{enumerate}
   \item $\Phi$ and $\Psi$ 
    commute with $z\mapsto z+2\pi i$;
   \item 
    $\re \Phi(z)\to \pm \infty$ as $\re z \to \pm \infty$  
      (and similarly for $\Psi$),
   \item for sufficiently large $R$,
          $\Phi(F^{-1}(\H_R))\subset \W$ and
          $\Phi^{-1}(G^{-1}(\H_R))\subset \V$, and \label{item:domains}
   \item
     $\Psi\circ F = G\circ \Phi$ wherever both compositions are defined.
  \end{enumerate}

 Let $\phi:\C\to\C$ be a quasiconformal map. Since
  $\phi$ is an order-preserving homeomorphism fixing $\infty$, we can
   define a branch of
   $\arg\phi(z) - \arg z$ in a full neighborhood
   of $\infty$. 
   It is well-known \cite[Lemma 4]{alexmisha}
   that there is some $C>1$ such that
   \begin{align}
     |z|^{1/C} \leq |\phi(z)|&\leq
         |z|^{C}\quad \label{eqn:hoelder}
    \text{and} \\
     |\arg \phi(z) - \arg z | 
                     &\leq C\log|z|  \label{eqn:complexhoelder}
   \end{align}
   when $z$ is sufficiently large.\footnote{%
   While there surely is a classical reference for 
    (\ref{eqn:complexhoelder}), 
    we were unable to locate one; Eremenko and Lyubich refer to
    \cite{lehtovirtanen}, but we did not find it there. 
    A short proof can be found in the Appendix of
    \cite{sebastianmisiurewiczspiral}.} 
    Translating this statement into
    logarithmic coordinates, we obtain the following fact.

 \begin{lem}[Hyperbolic distance of pullbacks]
    \label{lem:pullbacks} 
   Suppose that $F,G\in \Blog$ are normalized and quasiconformally equivalent.
    Then there are constants 
    $C>0$ and $M>0$ 
    such that 
    \[ \dist_{\H}(F_T^{-1}(z),G_{\wt{T}}^{-1}(w)) \leq
        C + \frac{\dist_{\H}(z,w)}{2} \]
    for all tracts $T$ of $F$ and $z,w\in\H_M$, 
    where $\wt{T}$ is the tract of $G$ containing $\Phi(F_T^{-1}(\H_M))$. 
 \end{lem}
 \begin{proof}[Sketch of proof]
  Let $\Psi,\Phi$ be the maps from the definition of quasiconformal
   equivalence.
   There are quasiconformal maps $\phi,\psi:\C\to\C$ such that
   $\phi\circ\exp = \exp\circ\Phi$ and
   $\psi\circ\exp = \exp\circ\Psi$. 
   Applying (\ref{eqn:hoelder}) and (\ref{eqn:complexhoelder})  
   to $\phi$ and $\psi^{-1}$, we easily see that there is some $M_0>0$ 
   such that
   $\dist_{\H}(z,\Phi(z))$ and
   $\dist_{\H}(z,\Psi^{-1}(z))$ are bounded,
   say by $\rho$,  when $z\in \H_{M_0}$. 

  By (\ref{eqn:hypexpansion}), we may also choose $M_1$ 
   sufficiently large so that
   $\|DF(z)\|_H\geq 2$ when $\re F(z)>M_1$. Finally, let
   $M>\max(M_0,M_1,R)$,
   where $R$ is as in part (\ref{item:domains}) of the
   definition of quasiconformal equivalence, be sufficiently large
   such that $\re z>M_0$ whenever $\Psi(F(z))\in \H_M$. 

  If $w\in\H_M$, we have
    $G_{\wt{T}}^{-1}(w) = \Phi(F_T^{-1}( \Psi^{-1}(w)))$, and hence
    \begin{align*}
      \dist_{\H}(F_T^{-1}(z),G_{\wt{T}}^{-1}(w)) &\leq
       \rho + \dist_{\H}(F_T^{-1}(z),F_T^{-1}(\Psi^{-1}(w))) \leq 
       \rho + \frac{\dist_{\H}(z,\Psi^{-1}(w))}{2} \\ 
     &\leq 
        \rho + \frac{\rho+\dist_{\H}(z,w)}{2} =
        3\rho/2 + \frac{\dist_{\H}(z,w)}{2} 
    \end{align*}
   when $z,w\in \H_M$. \end{proof}

 \begin{rmk}[Functions with quasidisk tracts] \label{rmk:quasidisks}
  It is not always easy to check whether two given functions are
   quasiconformally equivalent. However, suppose that 
   $U$ and $\tilde{U}$ are quasidisks whose boundaries contain $\infty$.
   Let 
   $f:U\to W$ and $g:\tilde{U}\to W$ be universal covering maps
   (where again
    $W=\C\setminus\cl{D}$ for a bounded Jordan domain $D$) that extend 
   continuously to the boundary of 
   $U$ resp.\ $\tilde{U}$ in $\C$. 

  Then we can pick a conformal isomorphism $\phi:U\to\tilde{U}$ such
   that $g\circ\phi = f$. Because $U$ and $\tilde{U}$ are
   quasidisks, $\phi$ extends to a quasiconformal map $\phi:\C\to\C$
   (see \cite[Satz 8.3]{lehtovirtanen} or
    \cite[Section 4.9]{hubbardteichmueller}).

  Hence, if $f,g\in\B$ are such that $f^{-1}(\{|z|>R\})$ and
   $g^{-1}(\{|z|>R|\})$ are quasidisks for large $R$, then $f$ and $g$
   are quasiconformally equivalent near infinity. More generally, if
   $F:\V\to\H$ is a function in $\Blog$ such that $\exp(\V)$ is a 
   quasidisk, then $F$ is quasiconformally equivalent to any
   function $G\in\Blog$ with the same property. 
   The tracts of the functions in Figures
   \ref{fig:julia} and \ref{fig:JR} are all quasidisks.
 \end{rmk}

\subsection*{External Addresses}
  Let $F\in\Blog$. We say that
   $z, w \in J(F)$
   \emph{have the same external address}
   (under $F$) if, for every $n\geq 0$, the points
   $F^n(z)$ and $F^n(w)$ belong to the closure of the same
   tract $T_n$ of $F$. 
 
  The sequence $\s= T_0 T_1 T_2 \dots$ is called the
   \emph{external address} of $z$ (and $w$) under
   $F$; compare \cite{strahlen} for a more detailed discussion.

   \begin{lem}[Expansion along orbits] \label{lem:expansion}
    Suppose that $F\in\Blog$ is normalized. If
     $z$ and $w$ have the same external address under
     $F$, then
    \[ |F^n(z) - F^n(w)| \geq 2^n|z-w| \]
     for all $n\geq 0$.
   \end{lem}
   \begin{proof} This is a direct consequence of the expansion property
     (\ref{eqn:expansion}). \end{proof}

 \subsection*{Further properties of quasiconformal maps}
  Throughout the article, we require a number of 
   well-known properties
   of quasiconformal maps. We collect a few of these here for the
   reader's convenience. By convention,
   the ``dilatation'' of a quasiconformal map $\psi$ will always mean
   the \emph{complex dilatation}; that is,
     \[ \dil(\psi) = \frac{\bar{\partial}\psi}{\partial\psi}. \]

 \begin{prop}[Compactness of qc mappings 
           {\cite[Sections II.5 and IV.5]{lehtovirtanen}}]
       \label{prop:compactness}
   Consider a sequence $\Psi_n:\C\to\C$ of quasiconformal maps,
   and suppose that there is a dense set $E\subset\C$ such that
   $(\Psi_n)$ stabilizes on $E$; i.e., for all $z\in E$ there is $n_0$
   such that $\Psi_n(z) = \Psi_{n_0}(z)$ for all $n\geq n_0$.

  If the maximal dilatation of the maps $\Psi_n$ is bounded
   independently of $n$, 
   then the sequence $\Psi_n$ converges locally uniformly to
   a quasiconformal map $\Theta:\C\to\C$. 

  If furthermore the complex dilatations $\dil(\Psi_n)$ converge
   pointwise almost everywhere, then their limit agrees with
   $\dil(\Theta)$ almost everywhere.
 \end{prop}

 \begin{prop}[Royden's Glueing Lemma {\cite[Lemma 2]{bersglueing}, \cite[Lemma 2]{polylikemaps}}] 
            \label{prop:glueing}
   Suppose that $U\subset\C$ is open,
      and that $\phi:U\to \phi(U)\subset\C$ is quasiconformal. Suppose
      furthermore that $\psi:\C\to\C$ is a quasiconformal map such that
      the function
           \[ \theta:\C\to\C;  z \mapsto \begin{cases}
                          \phi(z) \quad \text{if }z\in U; \\
                          \psi(z) \quad \text{otherwise} \end{cases} \]
      is a homeomorphism. Then $\theta$ 
         is quasiconformal. 
 \end{prop}


 \begin{prop}[QC maps of an annulus \cite{lehtoextension}] 
         \label{prop:interpolation}
   Let $A,B\subset\C$ be bounded annuli, 
      each bounded by two Jordan curves. 
      Suppose that $\psi,\phi:\C\to\C$ are quasiconformal maps such that
      $\psi$ maps the inner boundary $\alpha^-$ of $A$ to the inner boundary
      $\beta^-$ 
      of $B$, and
      $\phi$ takes the outer boundary $\alpha^+$ of $A$ to the
      outer boundary $\beta^+$ of $B$. 

     Let $z\in \alpha^-$ and $w\in \alpha^+$, let $\gamma$ be a curve in $A$
      connecting $z$ and $w$, and let $\tilde{\gamma}$ be a curve connecting
      $\psi(z)$ and $\phi(w)$ in $B$.

     Then there is a quasiconformal map $\wt{\phi}:\C\to\C$ that agrees with
      $\psi$ on the bounded component of $\C\setminus A$ and
      with $\phi$ on the unbounded component of $\C\setminus A$, and such that
      $\wt{\phi}(\gamma)$ is homotopic to 
      $\tilde{\gamma}$ relative $\partial B$. 
  \end{prop}
  \begin{remark}
   The statement about the homotopy class is not made in
    \cite{lehtoextension}, but follows directly from the proof. 
    (Alternatively, $\wt{\phi}$ can be obtained from \emph{any} 
     quasiconformal map that interpolates $\psi$ and $\phi$
     by postcomposition with a suitable Dehn twist.)
  \end{remark}

 Let us also formulate the translation of the preceding result to logarithmic
  coordinates, since we frequently use it in this setting. 
 \begin{cor}[Interpolation of quasiconformal maps] \label{cor:interpolation}
  Suppose that $H$ and $H'$ are $2\pi i$-periodic, unbounded Jordan domains,
   both containing some right half plane, 
   with $\cl{H'}\subset H$. 

  Suppose that $\Psi,\Phi:\C\to\C$ are quasiconformal maps, commuting with
   translation by $2\pi i$, such that $\cl{\Phi(H')}\subset \Psi(H)$. 
   Then there is a quasiconformal map $\wt{\Phi}:\C\to\C$ that agrees with
   $\Psi$ on $\C\setminus H$, agrees with $\Phi$ on $H'$, and commutes
   with translation by $2\pi i$. 
 \end{cor} 

  Finally, we will use the ``$\lambda$-lemma'' for
   holomorphic motions, as developed in 
   \cite{mss} and improved in \cite{bersroyden}; 
   compare \cite[Section 5.2]{hubbardteichmueller}.

  \begin{prop}[$\lambda$-lemma {\cite[Theorem 1]{bersroyden}}]
    \label{prop:lambdalemma}
   Let $E\subset\C$ and $R>0$, and suppose that the functions 
      \[ h_{\lambda}:E\to \C, \quad \lambda\in\D_{R}(0) \]
    are injective, with $h_{0}=\id$, 
    and furthermore depend holomorphically on
    $\lambda$ for fixed $z\in E$. (Under these assumptions,
    we say that the $h_{\lambda}$ form a \emph{holomorphic motion} of 
    the set $E$.)

   Then each $h_{\lambda}$ extends to a quasiconformal self-map of the plane.
    The complex
    dilatation of this map is bounded by $|\lambda|/R$.
  \end{prop}
  \begin{remark}[Remark 1]
   \cite[Section 5.2]{hubbardteichmueller} even establishes
    the stronger fact, due to
    Slodkowski, that the extensions of the $h_{\lambda}$ can themselves
    be chosen to
    depend holomorphically on $\lambda$. 
  \end{remark}
  \begin{remark}[Remark 2]
   If each $h_{\lambda}$ commutes with translation by $2\pi i$, then
    the extension can also be chosen with this property. (Apply the
    above theorem to the holomorphic motion $g_{\lambda}$ of
    $\exp(E)\cup\{0\}$ defined by $g_{\lambda}(0)=0$ and
    $g_{\lambda}(\exp(z)) := \exp(h_{\lambda}(z))$.)
  \end{remark}

\section{Conjugacy near infinity}
  \label{sec:conjugacy}

  In this section, we prove Theorem \ref{thm:main1}. We begin
   by treating the special case where both maps are 
   of disjoint type. 
 
 \begin{thm}[Conjugacy between disjoint-type maps]
  \label{thm:disjointconjugacy}
  Suppose that two functions in $\Blog$,  
    \[ F:\V\to H\quad\text{and}\quad
       G:\Phi(\V)\to\Psi(H) \]
   are quasiconformally equivalent, 
    $\Psi\circ F = G\circ\Phi$. Suppose furthermore that
   $F$ and $G$ are of disjoint type; i.e.,
   $\cl{\V}\subset H$ and $\cl{\Phi(\V)}\subset\Psi(H)$.
 
   Then there is a quasiconformal map $\Theta:\C\to\C$ with the
   following properties:
   \begin{enumerate}
    \item $\Theta|_{\V}$ is isotopic to
     $\Phi|_{\V}$ relative $\partial \V$. 
    \item $\Theta$ is a conjugacy between $F$ and $G$; i.e.\ 
       $\Theta\circ F = G\circ\Theta$ on $\V$.
    \item $\dil(\Theta)=0$ almost everywhere on $J(F)$.
    \item $\Theta( z +2\pi i)=\Theta(z)+2\pi i$.
   \end{enumerate}
 \end{thm}
 \begin{proof} By Corollary
   \ref{cor:interpolation} (picking a $2\pi i$-invariant unbounded Jordan 
   domain $H'$ with $\V\subset H'$ and $\cl{\Phi(H')}\subset \Psi(H)$), we
   can find a quasiconformal map $\wt{\Phi}:\C\to\C$ such that
   $\wt{\Phi}$ agrees with $\Phi$ on $\V$ and with
   $\Psi$ on $\C\setminus H$ (and such that
   $\wt{\Phi}$ still commutes with addition by $2\pi i$). Since
   $\Phi$ and $\wt{\Phi}$ agree on the domain of definition of $F$,
   we clearly have $\Psi\circ F = G\circ \wt{\Phi}$. 

  In analogous manner, we can modify $\Psi$ to a quasiconformal
   map $\Psi_0:\C\to\C$ that
   is conformal on a neighborhood of $\cl{\V}$, agrees
   with $\Psi$ on $\C\setminus H$, and commutes with addition by $2\pi i$. 
   (Compare also the main result of \cite{lehtoextension}.) Note that
   we are not claiming that this modified map $\Psi_0$ will satisfy
   the same functional equation as $\Psi$. 

  By the Alexander trick,
   the isotopy class of a homeomorphism between two Jordan domains is
   determined by its boundary values
   (compare also \cite[Proposition 6.4.9]{hubbardteichmueller}). 
   Hence the maps 
   $\Psi|_H$, $\wt{\Phi}|_H$ and $\Psi_0|_H$ all 
   belong to a single isotopy class relative $\partial H$.

 We now define a sequence of maps $\Psi_n:\C\to\C$ inductively,
  starting with $\Psi_0$, by setting
      \[ \Psi_{n+1}|_T := G^{-1}_{\Phi(T)}\circ \Psi_n \circ F|_T \]
  for every tract $T$ of $F$, and 
      \[ \Psi_{n+1}|_{\C\setminus\V} := \wt{\Phi}|_{\C\setminus \V}. \]

 Clearly each $\Psi_n$ is a homeomorphism (recall that the components
  of $\V$ accumulate only at infinity by definition). 
  By the glueing lemma (Proposition \ref{prop:glueing}), it follows that
  each $\Psi_n$ is quasiconformal. Since $F$ and $G$ are holomorphic,
  the maximal dilatation of
  $\Psi_n$ depends only on that of $\Psi_0$ and $\wt{\Phi}$, and is hence
  bounded independent of $n$. 

 Furthermore,
  $\Psi_n|_H$ is isotopic to
  $\Psi|_H$ relative $\partial H$ for all $n$.
  This implies that 
  the maps $\Psi_{n+1}|_{\V}$ and $\Phi|_{\V}$ are isotopic
  relative $\partial \V$. 

  By construction, $\Psi_n\circ F = G\circ \Psi_{n+1}$, and
  $\Psi_n$ and $\Psi_{n+1}$ agree outside of the
  set $F^{-n}(H)$, so the sequence $\Psi_n$ stabilizes on the set 
   \[ \C \setminus \bigcap_{n=0}^{\infty} F^{-n}(H) = \C \setminus J(F). \]
   By Lemma \ref{lem:emptyinterior}, $\C\setminus J(f)$ is a dense subset
   of $\C$, and it follows from
   Proposition \ref{prop:compactness} that 
   the $\Psi_n$ converge to some quasiconformal map
   $\Theta:\C\to\C$ with $\Theta\circ F = G \circ \Theta$. 

   The dilatations of the maps 
    $\Psi_n$ stabilize on the set $\C\setminus J(F)$,
    but on the other hand each $\Psi_n$ is conformal on a neighborhood
    of $J(F)$, so that its complex dilatation is zero there. In particular,
    the dilatations converge pointwise, and it follows from the
    second part of Proposition \ref{prop:compactness} that 
    $\dil(\Theta)=0$ almost everywhere on $J(F)$. 
  
 Furthermore,
  $\Theta|_{\V}$ belongs to the isotopy class of 
  $\Phi|_{\V}$ relative $\partial \V$. Since each
  $\Psi_n$ has $\Psi_n(z+2\pi i) = \Psi_n(z) + 2\pi i$, the same is
  true of $\Theta$. 
 \end{proof}

  Now let 
    \[ F_{0}:\V\to\H\]
   be an arbitrary normalized 
   function in $\Blog$. 
   We consider the one-dimensional family
   \[ F_{\kappa}: (\V - \kappa) \to \H; \quad z \mapsto
      F_0(z+\kappa) \quad\quad (\kappa\in\C). \]
   Note that all maps $F_{\kappa}$ are normalized. We will now prove
   Theorem \ref{thm:main1} for this family, which
   implies the general statement
   when combined with Theorem
   \ref{thm:disjointconjugacy};
   see Corollary \ref{cor:conjugacy} below.

  For given $\kappa\in\C$, we define maps
   $\Theta_n = \Theta_n^{\kappa}$ by
   $\Theta_0(z) := z$ and
   \[ \Theta_{n+1}(z) := (F_0)^{-1}_T(\Theta_n(F_0(z))) - \kappa
                              \quad\text{(wherever defined)}, \]
   where $T$ is the tract of $F_0$ containing $z$. In other words,
   $\Theta_n$ is obtained by iterating forward $n$ times under $F_0$, and then
   taking the corresponding pullbacks under $F_{\kappa}$. 

  \begin{thm}[Convergence to a conjugacy] \label{thm:nondisjointconjugacy}
   Let $\kappa\in\C$, and let $Q>2|\kappa|+1$. Then 
    the functions $\Theta_n$ are defined and continuous on 
    $J_Q(F_0)$, where they converge uniformly to
    a map
    \[ \Theta = \Theta^{\kappa} : J_Q(F_0) \to J(F_{\kappa}) \]
    that satisfies $\Theta\circ F_0 = F_{\kappa}\circ \Theta$,
      \begin{equation}
         |\Theta(z)-z|\leq 2|\kappa|  \label{eqn:distanceTheta}
      \end{equation}
    and is
    a homeomorphism onto its image. 

  For fixed $Q>1$ and $z\in J_Q(F_0)$, the 
    function $\kappa\mapsto \Theta^{\kappa}(z)$ is holomorphic
    on $\D_{(Q-1)/2}$.
  \end{thm}
 \begin{proof} The functions $\Theta_n$ are clearly continuous where defined. 
   Let us show inductively that $\Theta_n(z)$ is defined and
    \begin{equation}
     |\Theta_n(z) - z|\leq 2|\kappa| \label{eqn:distance}
    \end{equation}
   whenever $z\in J_Q(F_0)$. Indeed,
   for such $z$ we have  
   $\re F_0(z)\geq Q > 2|\kappa|+1$, so the induction hypothesis implies that
   $\Theta_n(F_0(z))\in\H$, and thus $\Theta_{n+1}(z)$ is defined. 
   Furthermore, by the expansion property 
    (\ref{eqn:expansion}) of $F_0$ and the induction hypothesis, we see that
    \begin{align*}
      |\Theta_{n+1}(z) - z| &=
      |(F_0)^{-1}_T(\Theta_n(F_0(z))) - \kappa -
       (F_0)^{-1}_T(F_0(z))| \\
       &\leq \frac{1}{2}|\Theta_n(F_0(z)) - F_0(z)| +
               |\kappa| \leq |\kappa| + |\kappa| = 2|\kappa|,
    \end{align*}   
    as required. 

  Using (\ref{eqn:distance}), we see that
    \begin{align*}
     |\Theta_{n+k}(z) &- \Theta_n(z) | = 
       |(F_0)_T^{-1}(\Theta_{n-1+k}(F_0(z))) - 
                    (F_0)_T^{-1}(\Theta_{n-1}(F_0(z)))| \\
       &\leq \frac{1}{2} |\Theta_{n-1+k}(z) - \Theta_{n-1}(F_0(z))|
       \leq \dots  \leq
       \frac{1}{2^n} |\Theta_k(F_0^n(z)) - \Theta_0(F_0^n(z)) | \leq
         \frac{2|\kappa|}{2^n}.
    \end{align*}
    Hence the functions $\Theta_n$ form a Cauchy sequence, and thus converge
    to some function
     \[ \Theta=\Theta^{\kappa}:J_Q(F_0) \to J_1(F_{\kappa}) \]
    satisfying (\ref{eqn:distanceTheta}) and
    $\Theta\circ F_0 = F_{\kappa}\circ \Theta$. 
    Since the convergence is locally  uniform in $\kappa$ and
     each $\Theta_n$ is holomorphic in $\kappa$, the map
     $\Theta$ likewise depends holomorphically on $\kappa$.

    It remains to verify that $\Theta$ has the stated properties. Note that,
     by definition of $\Theta$, the external address $\tilde{\s}$
     of $\Theta(z)$ under $F_{\kappa}$
     is
     determined uniquely by the external address $\s$ of $z$ under $F_0$.
     Indeed, if $\s=T_1 T_2 \dots$, then
     $\tilde{\s} = \tilde{T}_1 \tilde{T}_2 \dots$, where
     $\tilde{T}_j = T_j - \kappa$. 

     To see that $\Theta$ is injective, suppose that
     $\Theta(z)=\Theta(w)$. Then $z$ and $w$
     have the same external address under
     $F_0$, 
     and by (\ref{eqn:distanceTheta}) their orbits  are never separated by more
     than $4|\kappa|$. By Lemma \ref{lem:expansion}, this is impossible unless
     $z=w$; so $\Theta$ is indeed injective. 

    Furthermore,
     $\lim_{z\to\infty} \Theta(z)=\infty$, again by
     (\ref{eqn:distanceTheta}), so $\Theta$ extends to a continuous
     injective map on the compact space $J_Q(F_0)\cup\{\infty\}$, and thus
     is a homeomorphism onto its image.
 \end{proof}

\begin{lem}[Image of $\Theta$]   \label{lem:imageoftheta}
  Let $\kappa\in\C$, and let $Q$ and $\Theta$ be
   as in the preceding theorem.
   Then 
   $\Theta(J_Q(F_0))\supset J_{2Q}(F_{\kappa})$. 
\end{lem}
\begin{proof}
  Set $G_0 := F_{\kappa}$ and consider the family
   $G_{\lambda}(z) := G_0(z+\lambda)$; then $F_0 = G_{-\kappa}$. 
   Applying Theorem \ref{thm:nondisjointconjugacy} to this family,
    we obtain a map $\Theta':J_Q(F_{\kappa})\to J(F_0)$ satisfying
    $\Theta'\circ F_{\kappa} = F_0\circ \Theta'$
   and (\ref{eqn:distanceTheta}). Now, if $w\in J_{2Q}(F_{\kappa})$, then
     $z := \Theta'(w)$ satisfies
      \[ \re F^k(z) \geq \re F^k(w) - 2|\kappa| \geq 2Q - 2|\kappa| > Q. \]
     So $z\in J_Q(F_0)$. 
   The points $w$ and $w' := \Theta(z)$ have the same external address under
     $F_{\kappa}$. Furthermore, $F_{\kappa}^k(w')=\Theta(F_0^k(z))$ and
     $F_0^k(z) = \Theta'(F_{\kappa}^k(w))$ for all $k$, and hence
     \[ |F_{\kappa}^k(w) - F_{\kappa}^k(w')| \leq 
           |F_{\kappa}^k(w) - 
            \Theta'(F_{\kappa}^k(w)) | +
             |F_0^k(z) - \Theta(F_0^k(z)) | \leq 4K. \]
    So by Lemma \ref{lem:expansion}, 
     we have $w=w' = \Theta(z)\in \Theta(J_Q(F_0))$, as required.  
\end{proof}

 \begin{thm}[Quasiconformal extension and dilatation of $\Theta$] 
          \label{thm:qcextension}
  Let $\kappa\in\C$, and let $Q$ and $\Theta$ be
   as in Theorem
    \ref{thm:nondisjointconjugacy}. Then $\Theta$ extends
    to a quasiconformal map
    $\Theta:\C\to\C$. 
    This extension can be chosen such that
    $\Theta(z + 2\pi i)=\Theta(z)+2\pi i$, and such that 
    $\Theta|_{\V}$ is isotopic to
    $\Phi(z) := z - \kappa$ relative
    $\partial \V$. 

  Furthermore, the maximal dilatation of $\Theta$ on $J_{Q'}(F_0)$ tends to
   zero as $Q'\to\infty$. In particular, the dilatation of
   $\Theta$ is zero almost everywhere on $I(F_0)\cap J_Q(F_0)$.
 \end{thm}
 \begin{proof}
   The functions $\Theta=\Theta^{\kappa}$ define a holomorphic
    motion of the set $J_{Q}$. By the 
     $\lambda$-lemma (Proposition \ref{prop:lambdalemma}),
     each of
     these functions extends to a quasiconformal self-map $\Theta^{\kappa}$ 
     of the plane. 

   For abbreviation, let us set $J_Q^{\kappa} := J_Q(F_{\kappa})$,
    and also write $J_Q := J_Q^0$. 
    As pointed out in Remark 2 after Proposition
     \ref{prop:lambdalemma}, $\Theta$ can be chosen to commute with
     translation by $2\pi i$. Also, by (\ref{eqn:distanceTheta}),
         \[ \Theta\bigl( F_0(J_Q)\bigr) \subset
                \cl{\H_1}, \]
     so we can use Corollary \ref{cor:interpolation} to obtain
     a quasiconformal map $\Theta':\C\to\C$ that agrees with
     $\Theta$ on $F_0(J_Q)$, but is the identity on
     $\C\setminus\H$ (and is hence isotopic to the identity
     relative $\partial\H$). Consider the map $\Theta''$, defined by 
        \[ \Theta''(z) := (F_0)_T^{-1}(\Theta'(F_0(z))) - \kappa \]
      when $z$ belongs to a tract $T$ of $F$, and 
      $\Theta''(z)=\Phi(z)$ otherwise.
      This map is quasiconformal, isotopic to $\Phi$ relative
      $\partial \V$, and agrees with $\Theta'$, and  hence $\Theta$,
      on $J_Q(F_0)$. 

     To discuss dilatation, recall from Theorem \ref{thm:nondisjointconjugacy}
      that the maps 
         $\Theta^{\kappa}|_{J_{Q'}}$,
      for $Q'>Q$, define a holomorphic motion over
      the disk $\D_{(Q'-1)/2}(0)$ in $\kappa$-space.
     It follows from 
      the dilatation statement in the $\lambda$-lemma that 
      $\Theta|_{J_{Q'}}$ extends to a quasiconformal map
      with dilatation bounded by $2|\kappa|/(Q'-1)$. In particular,
     \[ \dil(\Theta) \leq 2|\kappa|/(Q'-1) \quad\text{a.\ e.\ on $J_{Q'}(F_0)$}; \]
    clearly this bound tends to $0$ as $Q'\to\infty$, as claimed.    
   
  Finally, recall that we have
      \[ \Theta \circ F_0^n = F_{\kappa}^n \circ \Theta \]
   on $J_Q$. Since $F_0$ and $F_{\kappa}$ are holomorphic, 
    we see that (for $Q'\geq Q$) the maximal dilatation of
    $\Theta$ on 
    \[ X^n_{Q'} :=
           \{ z\in J_Q:\, F_0^n(z)\in J_{Q'} \} \]
    is the same as the maximal dilatation of $\Theta$ on
    $J_{Q'}$, which tends to $0$ as $Q'\to\infty$. Since the
    bound is independent of $n$, the same is true for 
     \[ X_{Q'} := \bigcup_{n=0}^{\infty} X^n_{Q'}. \] 
    But
     $I_Q(F_0) = \bigcap_{Q'\geq Q} X_{Q'}$,
    so the dilatation of
   $\Theta$ on 
   $I_Q(F_0)$ is zero, as required. 
 \end{proof}

 We are now ready to prove Theorem \ref{thm:main1}, which we restate
  (with some additional details) in logarithmic coordinates.

 \begin{cor}[Conjugacy between qc equivalent maps] \label{cor:conjugacy}
  Suppose that $F,G\in\Blog$ are quasiconformally equivalent,
    $\Psi\circ F = G\circ\Phi$. For
   sufficiently large $Q>0$, there exists a quasiconformal map
   $\Theta$ such that
   \begin{enumerate}
     \item $\Theta|_{\V}$ is isotopic to $\Phi|_{\V}$ relative $\partial \V$.
     \item $\Theta\circ F = G\circ \Theta$ on $J_Q(F)$.
     \item $\Theta(J_Q(F))\supset J_{Q'}(G)$ for some $Q'>Q$.
     \item The dilatation of $\Theta$ is zero on $J_Q(F)\cap I(F)$. 
          \label{item:conjugacydilatation}
     \item $\Theta(z+2\pi i)=\Theta(z)+2\pi i$.
   \end{enumerate}
 \end{cor}
 \begin{proof} Let $\V$ and $\W$ be the domains of $F$ and $G$. By 
   restriction and conjugation, as discussed in Section
   \ref{sec:preliminaries}, we 
    may suppose without loss of generality that
   $F$ and $G$ are normalized, and that
   $\Phi(\V)\subset \W$. 

  Choose $K,L>0$ sufficiently large that
   \[
     F_0:\underset{=:\V_0}{\underbrace{\V+K}} \to \H; z\mapsto F(z-K)
      \quad\text{and}\quad
     G_0:\underset{=:\W_0}{\underbrace{\W+L}} \to\H;z\mapsto G(z-L)
   \]
   are of disjoint type, and that furthermore
    $\cl{\Phi(\V)}+L\,\subset \Psi(\H)$. 

   Now we can apply Theorem
    \ref{thm:disjointconjugacy} to
    obtain a quasiconformal conjugacy
    $\Theta_2$ between $F_0$ and $G_0$.

   Furthermore, we can apply Theorems \ref{thm:nondisjointconjugacy} and
    \ref{thm:qcextension}    
    to $F$ and $F_0$, as well as to $G$ and $G_0$, obtaining quasiconformal
    maps $\Theta_1$ and $\Theta_3$. It is easy to check that the function
    \[ \Theta := \Theta_3^{-1}\circ \Theta_2\circ \Theta_1 \] has the required
    properties. \end{proof}

\begin{proof}[Proof of Theorem \ref{thm:main1}]
  Suppose $f,g\in\B$ are quasiconformally equivalent near infinity, i.e.\ 
   \begin{equation}
    \psi(f(z)) = g(\phi(z))  \label{eqn:equivalenceagain} \end{equation}
  whenever $|f(z)|$ or $|g(\phi(z))|$ is large enough, with
    $\phi,\psi:\C\to\C$ quasiconformal. Without loss of generality,
    we may assume that $\phi(0)=0$ and $\psi(0)=0$ (otherwise
    we modify these maps inside some bounded disk, using Proposition
    \ref{prop:interpolation}). 

  Pick a logarithmic transform $F:\V\to H$, where we may assume that the
   disk $\exp(H)$ is chosen sufficiently large to ensure that
   (\ref{eqn:equivalenceagain}) holds 
   for $z\in \exp(\V)$. Let $\Phi:\C\to\C$ and $\Psi:\C\to\C$ be lifts
   of $\phi$ and $\psi$, respectively, under the exponential map.
   Then 
    \[ G := \Psi \circ F \circ \Phi^{-1} \]
    is a logarithmic transform of $g$, and $F$ and $G$ are quasiconformally
    equivalent by definition. (Note that we are not claiming that
    \emph{all} logarithmic transforms of $f$ and $g$ are quasiconformally
    equivalent.) We define $\theta$ by
      $\theta(\exp(z)) := \exp(\Theta(z))$, where $\Theta$ is the map
      from the previous theorem, and are done. 
\end{proof}

 We subdivided the proof of Theorem \ref{thm:main1} into two
  steps, using Theorem
  \ref{thm:disjointconjugacy} to reduce the problem 
  to the simpler family $F_{\kappa}$. 
  We remark that this would not be
  necessary if we were willing to forgo the statement that the dilatation
  of $\theta$ on the escaping set is zero. 

 Indeed, we can adapt the proof of Theorem \ref{thm:nondisjointconjugacy}
  to construct a suitable map $\Theta$ for 
  any two quasiconformally equivalent functions
  $F,G\in\Blog$. We sketch the argument in the following. 

 Letting $\Psi$ and $\Phi$ denote the maps from the definition of 
  quasiconformal equivalence, we set $\Theta_0(z) := z$ and define
  $\Theta_n$ inductively as follows. If $T$ is a tract of $F$ and
  $\wt{T}$ is the tract of $G$ that contains 
  $\Phi(F_T^{-1}(\H_M))$ for sufficiently large $M$, we define for $z\in T$:
     \[ \Theta_{n+1}(z) := G_{\wt{T}}^{-1}(\Theta_n(F(z))) \]
  (where defined).

 By virtue of Lemma \ref{lem:pullbacks}, the proof of Theorem
  \ref{thm:nondisjointconjugacy} goes through as before if we
  replace uniform convergence in the Euclidean metric by 
  uniform convergence in the hyperbolic metric. That is, for sufficiently
  large $Q$, the maps $\Theta_n$ are all defined on $J_Q(F)$ and converge
  uniformly to a map $\Theta:J_Q(F)\to J(G)$ that is a homeomorphism onto its
  image. 

 It is important to observe that,
  for fixed $F$, 
  the convergence is uniform not only in $z$ but also in $G$ 
  if $\Phi$ and $\Psi$ range over a compact
  set of quasiconformal mappings. Hence it follows that
  the conjugacy $\Theta$ still depends
  holomorphically on $G$ (which was not immediately clear from our original
  proof of Corollary \ref{cor:conjugacy}). We state this result
  formally for future reference.

 \begin{prop}[Analytic dependence of $\theta$]
  Let $f\in\B$. Let $M$ be a finite-dimensional complex manifold,
   together with a base point $\lambda_0\in M$. Suppose that 
   $(f_{\lambda})_{\lambda\in M}$ is a family of entire functions
   quasiconformally equivalent to $f$, with the equivalences given by 
     $\psi_{\lambda} \circ f = f_{\lambda} \circ \phi_{\lambda}$,
   where 
   $\psi_{\lambda_0}=\phi_{\lambda_0}=\id$, and 
   $\phi_{\lambda}$ and $\psi_{\lambda}$ depend analytically on 
   $\lambda$. 

  Let $N\ni \lambda_0$ be a compact subset of $M$. 
   Then there exists a constant $R>0$ 
   such that, for every $\lambda\in N$, there is an injective
   function
   $\theta=\theta^{\lambda}:J_R(f) \to J(f_{\lambda})$ with the following
   properties:
   \begin{enumerate}
    \item $\theta^{\lambda_0} = \id$,
    \item 
      $\theta^{\lambda}\circ f = f_{\lambda} \circ \theta^{\lambda}$ and
    \item  
     for fixed $z\in J_R(f)$, the function $\lambda\mapsto\theta^{\lambda}(z)$ 
     is analytic in $\lambda$ (on the interior of $N$). \qedoutsideproof
   \end{enumerate}
 \end{prop}

 In particular, we can again use the
  $\lambda$-lemma to show that $\theta^{\lambda}$ has a quasiconformal
  extension, as in Theorem \ref{thm:qcextension}.  If
  one was able furnish a direct proof of the statement 
   that the dilatation on the escaping set is zero ---
   our argument used the fact that the parameter space of the family 
   $F_{\kappa}$ is a parabolic surface, and hence does not generalize --- 
   then
   Theorem \ref{thm:disjointconjugacy} would no longer be required for the
   proof of Theorem \ref{thm:main1}.
 
 It is not difficult to show directly that the map $\Theta$ constructed
  above agrees with the map from Corollary \ref{cor:conjugacy}. 
  (In particular, it \emph{does} have zero dilatation on the escaping set.) 
  This also follows from the results proved in the next section 
  (see Corollary \ref{cor:uniqueness}).

 \section{Rigidity} \label{sec:rigidity}

 Let us now show that a (not necessarily quasiconformal) conjugacy
  between two quasiconformally equivalent maps $F,G\in\Blog$
  only moves escaping orbits 
  by a bounded hyperbolic distance, provided that it 
  ``preserves combinatorics'' (condition (\ref{item:sametract}) 
  below). This, together
  with the existence results from the previous section,
  will allow us to deduce a number of rigidity statements
  (Corollaries \ref{cor:uniqueness} and \ref{cor:linefields} and Theorems
   \ref{thm:nolinefields} and \ref{thm:qcrigidity}). 

 \begin{thm}[Restriction on conjugacies] \label{thm:conjugacyrestriction}
   Let $F,G\in\Blog$ be normalized and 
    quasiconformally
    equivalent, say $\Psi\circ F = G\circ \Phi$.
    Suppose that $Q>0$ and that 
    $\Pi:J_Q(F)\to J(G)$ is continuous such that
    \begin{enumerate}
     \item $\Pi\circ F = G\circ \Pi$,  \label{item:conjugacy}
     \item $\Pi(z)\to\infty$ as $z\to\infty$, 
     \item $\Pi(z+2\pi i)= \Pi(z)+2\pi i$, and  \label{item:Hcommutes}
     \item for every $z\in J_Q(F)$, 
            $\Pi(z)$ and $\Phi(z)$ belong to the
            same tract of $G$. \label{item:sametract}
    \end{enumerate}
   If $Q'$ is sufficiently large, then the hyperbolic distance
    $\dist_{\H}(z,\Pi(z))$ is uniformly bounded on
    $J_{Q'}(F)$. 
  \end{thm}
  \begin{remark}
   The hypothesis that $\Pi$ is defined on $J_Q(F)$ can be considerably weakened
    (with the same proof). For example, it would be sufficient
    to assume that $\Pi$ is defined and continuous on a forward
    invariant set $A\subset J_Q(F)$ with the property that
    $A$ contains the grand orbit (in $J_Q(F)$) of at least one 
    point $z_0$. 
  \end{remark}
  \begin{proof}
%
   Let $C,M>0$ be the constants
    from Lemma \ref{lem:pullbacks}; by enlarging $M$ if necessary
    we may assume that $M\geq Q$.
    By Corollary \ref{cor:unboundedescapingsets}, 
    we can choose some point $z_0\in J_Q(F)$ such that 
    $\re z_0 \geq M$ and $\re \Pi(z_0) \geq M$; we set
    \[ \rho := \max(2C, \dist_{\H}(z_0,\Pi(z_0))). \]    

  Set $Q' := e^{\rho}\cdot \re(z_0) + 2\pi > Q + 2\pi$.
  We will show that
    $\dist_{\H}(z, \Pi(z)) \leq \rho $
   for all $z\in J_{Q'}(F)$. 

 \begin{claim}
   For every
     $z\in J_{Q'}(F)$, there is a point $\zeta \in J_{Q}(F)$, belonging
     to the same tract of $F$ as $z$, with
     $|z - \zeta |<2\pi$ and $F(\zeta)\in \{z_0+2\pi i k:k\in\Z\}$. 
 \end{claim} 
 \begin{proof}[Proof of the claim]
  $F$ maps the boundary of the tract $T$ containing 
    $z$ to the imaginary axis, and the distance of $z$ to 
    $\partial T$ is at most $\pi$. Since
    $\re F(z) \geq Q' \geq \re z_0$, we can hence 
    find a point $\zeta_1\in T$ with
    $|z-\zeta_1|<\pi$ and      $\re(F(\zeta_1))=\re(z_0)$. 
     There is a point $\zeta_2\in \{z_0+2\pi i k:k\in\Z\}$ with
     $|F(\zeta_1) - \zeta_2|\leq\pi$. We set
     $\zeta := F_T^{-1}(\zeta_2)$. 
     By the expansion property (\ref{eqn:expansion})
     of $F$, we have $|\zeta - \zeta_1|\leq\pi/2$,
     and are done. \noqed
  \end{proof}

  Now let $z\in J_{Q'}(F)$. For each $n\geq 0$, we can apply the Claim
   to $F^n(z)$ to obtain a point $\zeta^n\in J_Q(F)$ with
   $|F^n(z)-\zeta^n|<2\pi$ and $F(\zeta^n)\in \{z_0+2\pi i k:k\in\Z\}$. 
   We now pull back $\zeta^n$ along the orbit of $z$ to obtain a point
   $z^n$; i.e.,
     \[ z^n = F_{T_0}^{-1}(F_{T_1}^{-1}(\dots F_{T_{n-1}}^{-1}(\zeta^n)\dots)), \]
     where $T_0 T_1 \dots$ is the external address of $z$. 
    By induction and the expansion property (\ref{eqn:expansion}), we have
       \begin{equation}
        |F^j(z) - F^j(z^n)| < \frac{2\pi}{2^{n-j}} 
                \label{eqn:distancebetweenzandzn}
       \end{equation}
     for $j=0,\dots,n$. In particular, $z^n\in J_Q(F)$ and $z^n\to z$. 

  We set
    $z^n_j := F^j(z^n)$ and 
    $w^n_j := \Pi(z^n_j) = G^j(\Pi(z^n))$. 
  Let us prove inductively that 
    \begin{equation}
     \dist_{\H}(z^n_j,w^n_j) \leq \rho    \label{eqn:distanceH}
    \end{equation}
    for $j=n+1,n,\dots,0$. 
   Indeed, we have $z^n_{n+1} = z_0 + 2\pi i k$ for some $k\in\Z$, and 
    hence 
    \[ \dist_{\H}(z^n_{n+1},w^n_{n+1}) =
       \dist_{\H}(z_0,\Pi(z_0)) \leq \rho, \]
    by property (\ref{item:Hcommutes}) and definition of $\rho$. 

   Furthermore, for $j\leq n$, we have
    \[ w^n_j = 
        G|_{\wt{T}}^{-1}(w_{j+1}^n) , \]
    where $\wt{T}$ is the tract of $G$ containing $w^n_j$. 
    By assumption (\ref{item:sametract}), 
     $\wt{T}$ is also the tract of $G$ containing
    $\Phi(z^n_j)$. 

   We observe that
    $z_{j+1}^n, w_{j+1}^n\in \H_M$. Indeed, if 
    $j=n$, this is true by choice of $z_0$. If $j<n$, 
    recall that $\re z_{j+1}^n\geq Q'-2\pi$ by
    (\ref{eqn:distancebetweenzandzn}) and
    $\dist_{\H}(z_{j+1}^n,w_{j+1}^n)\leq \rho$ by the induction
    hypothesis. Our choice of $Q'$ implies that
    $\re(w_{j+1})\geq \re z_0 \geq M$. 

   By Lemma \ref{lem:pullbacks} and the induction hypothesis, it
    follows that
    \[ \dist_{\H}(z^n_j,w^n_j) \leq
         C + \frac{\dist_{\H}(z^n_{j+1},w^n_{j+1})}{2} 
         \leq C + \frac{\rho}{2} \leq \rho, \]
   as claimed. 

  We have $z^n_0 = z^n \to z$, and hence by continuity of $\Pi$ also
   $w^n_0 = \Pi(z^n_0)\to \Pi(z)$. Therefore (\ref{eqn:distanceH})
   implies that $\dist_{\H}(z,\Pi(z))\leq \rho$, as desired. 
\end{proof}

 \begin{cor}[Uniqueness of conjugacies] \label{cor:uniqueness}
  Let $F$ and $G$ be quasiconformally equivalent. Then for every
   $Q>0$, there is $Q'\geq Q$ 
   with the following property. If
   $\Pi_1,\Pi_2:J_Q(F)\to J(G)$ are continuous functions satisfying
   the hypotheses (\ref{item:conjugacy}) to (\ref{item:sametract}) 
   of the previous
   theorem, then $\Pi_1(z) = \Pi_2(z)$ for all $z\in J_{Q'}(F)$. 
 \end{cor}
 \begin{proof} We may assume without loss of generality that
   $F$ and $G$ are both normalized. Let $Q'\geq Q$ be chosen
   such that $J_{Q'}(F)\subset \cl{I_{Q}(F)}$ (recall
   Corollary \ref{cor:unboundedescapingsets}). 

  It follows from Theorem
   \ref{thm:conjugacyrestriction} that there is $Q''\geq Q$ 
   such that, for all $z\in J_{Q''}(F)$,
   the points $\Pi_1(z)$ and $\Pi_2(z)$ have the same
   external address, and stay within bounded hyperbolic distance of
   each other. By the expansion property 
   (\ref{eqn:hypexpansion}) of $G$, this implies
   $\Pi_1(z)=\Pi_2(z)$. 

  So we have proved that $\Pi_1=\Pi_2$ on $J_{Q''}(F)$. Using 
   (\ref{item:sametract}),
   we see that $\Pi_1=\Pi_2$ on $I_Q(F)$. But
   $I_Q(F)$ is dense in $J_{Q'}(F)$, so we are done.
 \end{proof}

 \begin{cor}[No invariant line fields]  \label{cor:linefields}
  Let $F\in\Blog$. Then $F$ has no invariant line fields on its
   escaping set $I(F)$. 
 \end{cor}
 \begin{proof} Recall 
   that the existence of an invariant line field is equivalent to the
   existence of a non-zero $F$-invariant Beltrami form
   whose support is contained in $I(F)$ 
   \cite[Section 3.5]{mcmullenrenormalization}.

  So suppose that $\mu$ was such a Beltrami form. Recall that
    \[ I(F) = \bigcap_{Q>0} \bigcup_{n\geq 0} F^{-n}(J_{Q}(F)). \]
   Since $F$ is holomorphic, this implies that there is no $Q>0$ such that
   $\mu|_{J_Q(F)}$ is zero almost everywhere. Also observe that
   $2\pi i$-periodicity of $F$ implies that $\mu$ is $2\pi i$-periodic. 

  By the Measurable Riemann Mapping Theorem 
         \cite[Theorem 4.6.1]{hubbardteichmueller}, 
   $\mu$ gives rise to a quasiconformal homeomorphism $\Phi:\C\to\C$,
   which we may choose to commute with translation by $2\pi i$. 
   The map
     \[ G := \Phi\circ F \circ \Phi^{-1} \]
   is holomorphic, and clearly quasiconformally equivalent to $F$. 

  By Corollary \ref{cor:conjugacy}, there is a quasiconformal map $\Theta$,
   isotopic to $\Phi$ relative the boundary of the domain of definition
   $\V$ of $F$, which conjugates $F$ and $G$ on $J_Q(F)$, where $Q>0$ is
   sufficiently large. 

  By Corollary \ref{cor:uniqueness}, we then have
     \[ \Theta|_{J_{Q'}(F)} = \Phi_{J_{Q'}(F)} \]
    for sufficiently large $Q'$. Hence the dilatation of $\Theta$
    and $\Phi$ agree almost everywhere on $I_{Q'}(F)$. This is
    a contradiction: the dilatation of $\Theta$ on $I_{Q'}(F)$ is zero
    almost everywhere, but this is false for the dilatation
    $\mu$ of the map $\Phi$. 
  \end{proof}

 \begin{proof}[Proof of Theorem \ref{thm:nolinefields}]
  Let $f\in\B$, and let $F$ be a logarithmic transform of $f$. 
   If $f$ supported an invariant line field on its escaping set,
   then the same would be true for $F$. (As in the proof of 
   Corollary \ref{cor:linefields}, the support of the line field has
   nontrivial intersection with every set of the form
    $\{ z\in I(f) : |f^n(z)|\geq R \}$, 
    $R>0$.) Hence the theorem follows from Corollary \ref{cor:linefields}. 
 \end{proof}

 \begin{proof}[Proof of Theorem \ref{thm:qcrigidity}]
  Suppose that $f$ and $g$ are entire function with finitely many
   singular values, let $\pi:\C\to\C$ be a topological
   conjugacy betweeen $f$ and $g$, and let $\mathcal{O}$ be the orbit
   of some escaping point $z_0\in I(f)$. 

  For simplicity, let us assume that $f(0)=0$, and that $\pi(0)=0$. 
   This is no loss of generality, since
   any $f\in\B$ has infinitely many fixed points
   (see 
     \cite[Lemma 69]{adamthesis} or \cite{alexmisharussian4}; compare also
     \cite{langleyzheng} for a more general result).  
   However, we would like to point out
   that this assumption is not essential for the proof, and
   is made purely for convenience. 

  Let $S := S(f)\cup\{0\}$. We can pick a quasiconformal homeomorphism
   (in fact, a diffeomorphism) $\psi:\C\to\C$ that is isotopic to
   $\pi$ relative $S$. Using the functional relation
   $\pi\circ f = g\circ \pi$, the isotopy between $\pi$ and $\psi$ lifts 
   to an isotopy between $\pi$ and and a quasiconformal map $\phi:\C\to\C$ with
      \[ \psi \circ f = g \circ \phi. \]
    (Compare also \cite[Section 3]{alexmisha}.) In particular, 
    $f$ and $g$ are quasiconformally equivalent. 

  Now, as usual, we change to logarithmic coordinates: we let
   $F:\V\to H$ be a logarithmic transform of $F$, and $\Pi$ be a lift
   of $\pi$; i.e., $\pi\circ\exp = \exp\circ\Pi$. Then
     $G := \Pi\circ F \circ \Pi^{-1}$ is a logarithmic transform of $g$.

  The isotopies between $\pi$ and $\psi$  resp.\ $\phi$ lift to 
   isotopies between $\Pi$ and maps $\Psi$, $\Phi$ satisfying
     $\Psi \circ F = G\circ \Phi$,
    so $F$ and $G$ are quasiconformally equivalent as elements of $\Blog$. 

  Furthermore, if $M>0$ is sufficiently large, then no point
   $z\in \H_M$ leaves the domain $H$ under the isotopy between
   $\Pi$ and $\Psi$. It follows that, if $T$ is a tract of $F$ and
   $z\in T$ with $F(z)\in \H_M$, then 
      $\Phi(z) \in \Pi(T)$. 

   Let $\Theta$ be the map from Corollary \ref{cor:conjugacy}. Then by 
    Corollary
    \ref{cor:uniqueness}, we have   
      \[ \Theta|_{J_{Q'}(F)} = \Pi|_{J_{Q'}(F)} \]
    for some ${Q'}\geq 0$. If $\theta$ is the quasiconformal map
    defined by $\theta\circ\exp = \exp\circ\Theta$, then 
    $\theta$ and $\pi$ agree on the set
     \[ J_{e^{Q'}}(f) = 
          \{ z\in\C: |f^n(z)|\geq e^{Q'}\text{ for all $n\geq 1$} \}. \]

   Pick $k_0\in\N$ such that $f^{k_0}(z_0)\in J_{e^{Q'}}(f)$. Then 
   $\pi$ agrees with the quasiconformal map $\theta$ on the tail 
    $\mathcal{O}_{k_0}:= \{f^k(z_0):k\geq k_0\}$ of the orbit $\mathcal{O}$. 

  We can modify 
   the map $\theta$ (e.g.\ using Proposition \ref{prop:interpolation})
   on a compact subset of $\C\setminus \mathcal{O}_{k_0}$ to a quasiconformal
   function that maps $f^k(z_0)$ to $\pi(f^k(z_0))$ for $0\leq k < k_0$.
   This is the desired quasiconformal extension of $\pi|_{\mathcal{O}}$. 
 \end{proof}
 \begin{remark}
  Note that the assumption that $S(f)$ is finite was used only 
   to find a quasiconformal map $\psi$ isotopic to $\pi$. Hence we can
   weaken the assumptions of Theorem \ref{thm:qcrigidity} to require only that
   $f,g\in\B$ and that the conjugacy $\pi$ is isotopic, relative $S(f)$,
   to a quasiconformal
   self-map of the plane.
 \end{remark}
   
\section{Hyperbolic Maps} \label{sec:hyperbolic}
 
 Recall that $f\in\B$ is \emph{hyperbolic} if $S(f)$ is contained
  in the union of attracting basins of $f$. Since $S(f)$ is
  compact by definition, there are then only finitely 
  many such basins, which together make up the
  Fatou set. In particular, $f$ is hyperbolic if and only if the
  postsingular set
   \[ \P(f) = \cl{\bigcup_{j\geq 0} f^j(S(f))} \]
  is a compact subset of the Fatou set. 

In the following, we assume without loss of generality
  that $0$ is one of the attracting periodic points of $f$.

 We will show that such $f$ is semi-conjugate on its Julia set to
  a disjoint-type map quasiconformally equivalent to $f$, 
  and this semi-conjugacy is
  a conjugacy when restricted to the escaping set. In view of
  Theorem \ref{thm:disjointconjugacy}, this implies that
  any two
  hyperbolic maps that are quasiconformally equivalent near infinity
  are in fact topologically conjugate on their sets of escaping points,
  and hence proves Theorem \ref{thm:hyperbolicconjugacy}.

 It is easy to see that there is a bounded open neighborhood $U$ of the
  postsingular set $\P(f)$ such that $\cl{f(U)}\subset  U$. We set
  $W := \C\setminus \cl{U}$ and $V:=f^{-1}(W)\subset W$. Then
   \[ f:V\to W \]
  is a covering map, and hence expands the hyperbolic metric of
  $W$. We claim that this map is in fact \emph{uniformly} expanding. 
  (Compare also \cite[Theorem C]{ripponstallardhyperbolic}.)

 \begin{lem}[Uniform expansion] \label{lem:hyperbolicity}
  There is $C>1$ such that
    $\|Df(z)\|_{W} \geq C$
   for all $z\in V$. 
 \end{lem}
 \begin{proof} Since $f$ is a covering map, we just 
   need to show that the inclusion
   $i:V\to W$ is uniformly contracting. Since the density of the
   hyperbolic metric of $V$ tends to $\infty$ near $\partial V$, and $V$
   and $W$ have no common finite boundary points, it is sufficient
   to prove that 
   $\rho_W(z)/\rho_V(z)\to 0$ as $z\to \infty$. 

  The hyperbolic density of 
   $W$ satisfies $\rho_W(z)= O(1/(|z|\log |z|))$ as $z\to\infty$. 
   We now estimate the hyperbolic 
   metric of $V$, using Lemma \ref{lem:beardonpommerenke}. Fix some
   point $w\in \C\setminus W = \cl{U}$ such that $w$ belongs to the
   unbounded component of $\C\setminus S(f)$. 

 \begin{claim} There
   is a constant $C$ and a 
   sequence $(w_j)$ of (pairwise distinct) preimages of
   $w$ under $f$ such that $|w_{j+1}| \leq C|w_j|$ for
   all $j$.
 \end{claim}
 \begin{proof}[Proof of the claim]
  Pick a Jordan curve
   $\gamma$ surrounding $S(f)$, but not surrounding $w$, 
   and let $G$ be the unbounded component of $\C\setminus \gamma$. 
   If $\wt{G}$ is a component of $f^{-1}(G)$, then 
   $f:\wt{G} \to G$ 
   is a universal covering. Hence we can find a
   sequence $(w_j)$ of preimages of $w$ in $\wt{G}$ such that
   the hyperbolic distance 
   (in $\wt{G}$) between $w_j$ and $w_{j+1}$ is constant.
   By the standard estimate (\ref{eqn:standardestimate}) 
   on the hyperbolic distance in the
   simply connected domain $\wt{G}$, this implies that
   $|w_{j+1}|\leq C|w_j|$
   for some $C$ and sufficiently large $j$, as desired.  \noqed
 \end{proof}

 By Lemma \ref{lem:beardonpommerenke}, the hyperbolic
  metric of the domain $V' := \C\setminus\{w_n:n\in\N\}$ satisfies
  $1/\rho_{V'}(z) = O(|z|)$. Since
  $\rho_V \leq \rho_{V'}$ by Pick's theorem,
  this means that $\rho_W(z)/\rho_V(z)\to 0$ as
  $z\to\infty$, as claimed. \end{proof}

 Let $K\geq 1$; if $K$ is chosen
 sufficiently large, then  $\cl{U}\subset\D_{K/2}(0)$. 
 Furthermore, choose
 $R\geq K$ such that  
   \[ f^{-1}(\{|z|>R\})\subset \{|z|>K+1\}. \]
 We define $M := R/K$ and
  $g(z) := f(z/M)$. Then $g$ is of disjoint type. Indeed, 
  we have 
   $\U := g^{-1}(\{|z|>R\})\subset \{|z|>R+M\}$. We define
 \begin{align*}
 &\V := f^{-1}(\{|z|>R\}), \quad \wt{\V} := f^{-1}(\{|z|>K\})\quad\quad
               \text{and} \\
 \U_j := g^{-j}(\{|z|>R\}), \quad &\V_j := f^{-j}(\{|z|>R\}), \quad
    \wt{\V}_j := f^{-j}(\{|z|>K\}). \end{align*}
  Note that $\V_j\subset \wt{\V}_j\subset W$ for all $j$. 

 We now define a sequence $\theta_k$, where
   $\theta_0=\id$ and 
     \[ \theta_k : \U_{k-1}\to \wt{\V}_{k-1} \]
  is a conformal isomorphism for $k\geq 1$, such that
   \[ f(\theta_{k+1}(z)) = \theta_{k}(g(z)). \]

  Begin by setting $\theta_1(z) := z/M$. Furthermore, for $z\in \U_0$, let
   $\gamma_1(z)\subset \wt{\V}_0$ be the straight line segment connecting 
   $z = \theta_0(z)$ and $z/M = \theta_1(z)$. 

  To define $\theta_2$ let $z\in \U_1$. Since
   \[ f(\theta_1(z)) = \theta_0(g(z)), \]
   the curve $\gamma_1(g(z))$ has a preimage component 
   $\gamma_2(z)\subset \wt{\V}_1$ 
   under $f$ with one endpoint at $\theta_1(z)$; we
   define $\theta_2(z)$ to
   be the other endpoint. Then $f(\theta_2(z))=\theta_1(g(z))$. 

 We continue inductively: the curve
  $\gamma_{j+1}(z)\subset \wt{\V}_j$ 
  is the pullback of $\gamma_j(g(z))$ with
  one endpoint at $\theta_j(z)$, and $\theta_{j+1}(z)$ is defined as
  the other endpoint of this curve.

 It follows from the definition that each $\theta_{k+1}$ is continuous.
  Hence, for every component $G$ of $\U_{k}$, $(\theta_{k+1})|_G$ is 
  a branch of $f^{-1}\circ \theta_k\circ g$, and hence a conformal isomorphism
  onto some component of $\wt{\V}_k$. It is likewise easy to check that
  $\theta_{k+1}$ is surjective, so these maps are indeed conformal isomorphisms
  between $\U_{k}$ and $\wt{\V}_k$. 

 We furthermore note that $\theta_k(\U_k) = \V_k$ by the inductive 
  construction.

 \begin{thm}[Convergence to a semiconjugacy] \label{thm:hyperbolicsemiconjugacy}
  In the hyperbolic metric of $W$, the maps 
   $\theta_k|_{J(g)}$ converge uniformly to a continuous surjection
  \[ \theta: J(g)\to J(f) \]
  with $f\circ\theta = \theta\circ g$ and
   $\theta(I(g))=\theta(I(f))$. Furthermore, 
   $\theta:I(g)\to I(f)$ is a homeomorphism. 
 \end{thm}
 \begin{proof} Let $z\in \U_{k}$. By definition,
   \[ \dist_{W}(\theta_{k+1}(z),\theta_{k}(z)) \leq
      \ell_W(\gamma_{k+1}(z)). \]
    We have
    \begin{align*} \ell_W(\gamma_1(z)) &\leq
       \ell_{\{|w|>K/2\}}(\gamma_1(z)) \\ &=
      \log(1 + \log M / \log(2|z|/MK)) \leq \log(1+\log M/\log 2) =: \mu 
    \end{align*}
    for all $z\in \U_0$. 
    Since $\gamma_{k+1}(z)$ is obtained from $\gamma_1(g^k(z))$ by
    a branch of $f^{-k}$, and $f$ is uniformly expanding on $W$
    by Lemma \ref{lem:hyperbolicity}, we see that
     \begin{equation}
       \dist_W(\theta_{k+1}(z),\theta_k(z)) \leq \mu/C^k
                \label{eqn:dist_limiting}
     \end{equation}
    for all $z\in \U_k$. 

    In particular,
    the maps $\theta_k|_{J(g)}$ form a Cauchy sequence, and by the
    completeness of the hyperbolic metric have a (continuous) limit 
    \[ \theta:J(g)\to W. \]
    By (\ref{eqn:dist_limiting}), $\theta$ satisfies
       \begin{equation} \label{eqn:dist_limit}
          \dist_W(\theta(z),z)\leq \mu\cdot \frac{C}{C-1}. 
       \end{equation}

    By definition, if $z\in J(g)$, then
     $f^k(\theta(z))=\theta(g^k(z))\in W$ for all $k\in\N$.
     Hence $\theta(z)\in J(f)$. 
     Also
     note that, by (\ref{eqn:dist_limit}),
      $\theta(z_n)\to\infty$ if and only if
     $z_n\to\infty$, so $\theta$ maps escaping points to
     escaping points. 

    The map $\theta:I(g)\to I(f)$ is clearly surjective. Indeed,
     if $w\in I(f)$, then $w\in \wt{V}_k$ for all sufficiently large $k$.
     Any limit point $z$ of the sequence $z_k = \theta_k^{-1}(w)$ will
     have $\theta(z) = w$. (Note that $(z_k)$ cannot diverge to infinity
     by (\ref{eqn:dist_limiting}).) 

   To prove injectivity on $I(g)$, suppose by contradiction that
    $\theta(z_1) = \theta(z_2)$, where $z_1,z_2\in I(g)$, $z_1\neq z_2$. 
    It follows from the construction of $\theta$ that
    also $g^j(z_1)\neq g^j(z_2)$ for all $j\geq 0$. However, 
    $\theta$ is injective on a set of the form
     \[ J_{R'}(g)\cap I(g) = 
      \{z\in I(g) : |g^j(z)|\geq R'\text{ for all $j\geq 1$} \}. \]
    (This follows from 
     Corollary \ref{cor:uniqueness}, or alternatively from an argument
     analogous to the proof of injectivity in  Theorem
     \ref{thm:nondisjointconjugacy}.) 
     Since $g^j(z_1)$ and $g^j(z_2)$ belong to $J_{R'}(g)\cap I(g)$ for
     sufficienly large $j$, we obtain the desired contradiction.
     The details are left to the reader.

    Finally, $\theta(J(g))\cup\{\infty\}$ is the continuous
     image of a compact set, and thus itself compact. Since
     $I(f)\subset \theta(J(g))\subset J(f)$ and
     $J(f)\subset\cl{I(f)}$ by \cite{alexescaping}, we see that 
     $\theta$ is surjective. Compactness of $J(g)\cup\{\infty\}$ 
     and the fact that $\theta^{-1}(I(f))=I(g)$ imply that the image
     of any relatively closed subset of $I(g)$ under $\theta$ is relatively
     closed in $I(f)$. Hence $(\theta|_{I(g)})^{-1}$ is continuous.
   \end{proof} 

   Recall that, by a ``pinched Cantor Bouquet'', 
    we mean a metric space that is
    the quotient of a
    \emph{straight brush} in the sense of \cite{aartsoversteegen} by
    a closed equivalence relation on its endpoints. 
    As a corollary of Theorem \ref{thm:hyperbolicconjugacy}, 
    we obtain the following. 

   \begin{cor}[Pinched Cantor Bouquets]
     Let $f\in\B$ be hyperbolic and of finite order. Then every
      dynamic ray of $f$ lands, and the Julia set is a
      pinched Cantor Bouquet.
   \end{cor}
   \begin{proof}
     Bara{\'n}ski \cite{baranskihyperbolic} proved that, in the
    disjoint case, the Julia set is a straight brush, where all points
    except (some of) the endpoints of the brush belong to $I(f)$. The Corollary
    then follows immediately from our Theorem 
    \ref{thm:hyperbolicsemiconjugacy}. 
   \end{proof}

\appendix        

\section{Structure of Escaping Sets} \label{sec:fatoueremenko}

 In this section, we discuss the bearing our results have on some
  intriguing questions about escaping sets of entire functions
  that go back to Fatou's original 1926 article \cite{fatou},
  and Eremenko's study of the escaping set \cite{alexescaping}.
  Fatou observed that the Julia sets of certain explicit entire functions
  contain curves on which the iterates tend to $\infty$, and asked whether
  this property holds for much more general functions. Eremenko showed
  that (for an arbitrary entire function $f$), every component of
  the closure $\cl{I(f)}$ is unbounded. He then asks whether, in fact,
  every component of $I(f)$ is unbounded, and also whether
  every point of $I(f)$ can be connected to $\infty$ by a 
  curve in $I(f)$. (For a more detailed discussion of these
  questions and their history, compare \cite{strahlen}.)

 These questions suggest the study of the following properties for 
  an entire function $f$. 
 \begin{enumerate}
  \item[(F)] \emph{Fatou property:} There is a 
    curve to $\infty$ in $I(f)$;
  \item[(E)] \emph{Eremenko property:} Every connected
    component of $I(f)$ is unbounded;
  \item[(S)] \emph{Strong Eremenko Property:} Every point 
    $z\in I(f)$ can be connected to $\infty$ by a curve in $I(f)$.
 \end{enumerate}

 It is shown in \cite{strahlen} that there exist hyperbolic
   functions
   $f\in\B$ for which the Julia set contains no curves to $\infty$. Thus 
   property (F) (and, in particular, property (S)) can fail for functions
   in class $\B$. In fact, there are even hyperbolic functions
   whose Julia set contains no nontrivial curves at all.
   Together with Theorem \ref{thm:main1}, 
   this implies the following.
  \begin{cor}[No curves in the escaping set] \label{cor:nocurves}
   There exists an entire function $f\in \B$ with the following property. 
   If $g\in\B$ is quasiconformally equivalent to $f$ near infinity, then
   the escaping set $I(g)$ contains no nontrivial curves.
  \end{cor}
 \begin{proof}
  Let $f$ be the example constructed in \cite[Theorem 8.4]{strahlen}, 
   whose Julia 
   set contains no nontrivial curve. If $g$ is quasiconformally
   equivalent to $f$ near infinity, then by Theorem \ref{thm:main1}, 
   for sufficiently large $R$ the set $J_R(f)$ of points whose forward
   orbits
   are contained in $\C\setminus\D_R(0)$ is homeomorphic to a subset of
   the Julia set of $f$. Therefore $J_R(f)$ contains no nontrivial curve
   either. Since the image of a nontrivial curve under $f$ is again a
   nontrivial curve, the same holds for all sets $f^{-n}(J_R(f))$,
   $n\geq 0$.

  Suppose, by contradiction, that $I(g)$ does contain a nontrivial curve
   $\gamma:[0,1]\to I(g)$; we may assume that $\gamma$ is not constant
   on any interval. For every $n$, $\gamma^{-1}(f^{-n}(J_R(f)))$ is
   a closed subset of $[0,1]$ that contains no intervals, and hence
   is nowhere dense. However, we have
    \[ [0,1] = \gamma^{-1}(I(f)) \subset \bigcup \gamma^{-1}(f^{-n}(J_r(f))), \]
   which contradicts the Baire category theorem. 
 \end{proof}

 In \cite{eremenkoproperty}, we establish
  Eremenko's property for 
  every hyperbolic function $f\in\B$, and more generally any function
  $f\in\B$ with bounded postsingular set. This shows that
  a situation as in Corollary 
  \ref{cor:nocurves} cannot occur for property (E).
  Whether
  there is any entire function for which property (E) fails remains
  an open problem. We remark that
  even in 
  the exponential family, the study of connected components of
  $I(f)$ is far from trivial: while in the hyperbolic case, each
  such component consists of a single dynamic ray
  \cite{accessible}, this is already false for postsingularly periodic
  (``Misiurewicz'') exponential maps \cite{misindecomposable}. 

 Note that our Theorem \ref{thm:hyperbolicconjugacy} also shows that 
  \begin{quotation}
   \emph{for any quasiconformal equivalence class in
    class $\B$, each of the properties (F), (E) and (S) either holds for
    all hyperbolic maps or fails for all hyperbolic maps.}
  \end{quotation}
 
 Now consider the following uniform variants of the above properties.
 \begin{enumerate}
  \item[(UE)] 
    For every $z\in I(f)$, there exists some unbounded connected set
    $A\ni z$ such that
    $f^n|_A\to\infty$ \emph{uniformly}.
  \item[(US)]
    Every $z\in I(f)$ can be connected to $\infty$ by
    a curve $\gamma$ such that 
    $f^n|_{\gamma}\to\infty$ \emph{uniformly}.
 \end{enumerate} 
  In many proofs of the Eremenko Property, or
   the Strong Eremenko Property,
   they are in fact established in this uniform sense.
   It is possible, following the construction in
   \cite{strahlen}, to construct an entire function for which property
   (UE) fails.  

  Theorem \ref{thm:main1} shows that
  \begin{quotation}
   \emph{for any quasiconformal equivalence class of
    Eremenko-Lyubich functions,
    each of the properties (UE) and (US) either holds for
    all maps or fails for all maps.}
  \end{quotation}
 
  In \cite{strahlen}, property (US)
   is
   established for a large subset of $\B$, in particular for those
   of finite order (as well as finite compositions of such
   functions). The above-mentioned recent results of Bara{\'n}ski
   \cite{baranskihyperbolic} also imply this property
   for \emph{disjoint-type} functions $f\in\B$ of finite order
   (i.e., hyperbolic maps with a single fixed attractor). Hence Theorem
   \ref{thm:main1}, together with \cite{baranskihyperbolic}, 
   provides an alternative
   proof of property (US) --- and thus a 
   positive answer to Fatou's and Eremenko's questions ---
   for functions $f\in\B$ of finite order.

{\small
\bibliographystyle{hamsplain}
\bibliography{../../Biblio/biblio}}

\end{document}